\documentclass[11pt,a4paper]{article}
\title{Families of G-Constellations over Resolutions of Quotient Singularities}
\author{Timothy Logvinenko}	
\date{}
\usepackage{amssymb}
\usepackage{amsmath}
\usepackage{amsthm}
\usepackage{graphicx}
\usepackage[all]{xy}

\DeclareMathOperator{\homm}{Hom}

\DeclareMathOperator{\autm}{Aut}

\DeclareMathOperator{\gl}{GL}
\DeclareMathOperator{\gsl}{SL}

\DeclareMathOperator{\picr}{Pic}

\DeclareMathOperator{\cl}{Cl}

\DeclareMathOperator{\spec}{Spec\;}

\DeclareMathOperator{\hilb}{Hilb\;}

\DeclareMathOperator{\ext}{Ext}

\DeclareMathOperator{\supp}{Supp}
\DeclareMathOperator{\ann}{Ann}
\DeclareMathOperator{\ev}{ev}
\begin{document}
\def\bv{\mathbf{v}}
\def\kgc_{K^*_G(\mathbb{C}^n)}
\def\kgchi_{K^*_\chi(\mathbb{C}^n)}
\def\kgcf_{K_G(\mathbb{C}^n)}
\def\kgchif_{K_\chi(\mathbb{C}^n)}
\def\gpic_{G\text{-}\picr}
\def\gcl_{G\text{-}\cl}
\def\trch_{{\chi_{0}}}
\def\regring{{R}}
\def\regrep{{V_{\text{reg}}}}
\def\givrep{{V_{\text{giv}}}}
\def\lbar{{(\mathbb{Z}^n)^\vee}}
\def\genpx_{{p_X}}
\def\genpy_{{p_Y}}
\def\genpcn_{p_{\mathbb{C}^n}}
\theoremstyle{definition}
\newtheorem{defn}{Definition}[section]
\newtheorem*{defn*}{Definition}
\newtheorem{exmpl}[defn]{Example}
\newtheorem*{exmpl*}{Example}
\newtheorem{exrc}[defn]{Exercise}
\newtheorem*{exrc*}{Exercise}
\newtheorem*{chk*}{Check}
\newtheorem*{remarks*}{Remarks}
\theoremstyle{plain}
\newtheorem{theorem}{Theorem}[section]
\newtheorem*{theorem*}{Theorem}
\newtheorem{prps}[defn]{Proposition}
\newtheorem*{prps*}{Proposition}
\newtheorem{cor}[defn]{Corollary}
\newtheorem*{cor*}{Corollary}
\newtheorem{lemma}[defn]{Lemma}
\newtheorem*{claim*}{Claim}
\numberwithin{equation}{section}
\maketitle

\begin{abstract}

 Let G be a finite subgroup of $\gl_n(\mathbb{C})$. A study is made of 
the ways in which resolutions of the quotient space $\mathbb{C}^n / G$ can 
parametrise $G$-constellations, that is, $G$-regular finite length sheaves. 
These generalise $G$-clusters, which are used in the McKay correspondence 
to construct resolutions of orbifold singularities.

  A complete classification theorem is achieved, in which all the
natural families of $G$-constellations are shown to correspond to
certain finite sets of $G$-Weil divisors, which are a special sort of
rational Weil divisor, introduced in this paper. Moreover, it is shown 
that the number of equivalence classes of such families is always finite. 

Explicit examples are computed throughout using toric geometry.

\end{abstract}

\setcounter{section}{-1} 

\section{Introduction} \label{section-intro}

  Let $G \subseteq \gsl_3(\mathbb{C})$ be a finite subgroup and let
$X$ be the quotient space $\mathbb{C}^3 / G$. Nakamura made a study of
$G$-clusters, the $G$-invariant subschemes of dimension $0$ whose
coordinate ring, with the induced $G$-action, is the regular
representation $\regrep$ of $G$. He introduced the scheme $G$-$\hilb$,
which parametrises all $G$-clusters and showed \cite{Nak00} that, in
the case of $G$ being abelian, it is a crepant resolution of
$\mathbb{C}^3 / G$, conjecturing that the same holds for the
non-abelian case. 

  Craw and Reid \cite{Craw02} introduced an alternative way of explicit
calculation of $G$-$\hilb \mathbb{C}^3$ and in his thesis
\cite{Craw-thesis} Craw introduced the concept of $G$-constellation as a
generalisation of $G$-cluster. A $G$-constellation is a
$G$-equivariant coherent sheaf whose global sections
form the regular representation of $G$. In particular, the structure
sheaf of any $G$-cluster is a $G$-constellation.  

  $G$-constellations can be interpreted in terms of representations of
the McKay quiver of $G$. This allows for the use of an earlier result of King 
\cite{King94} on GIT construction of moduli spaces of quiver 
representations to introduce the stability conditions known as 
$\theta$-stability on $G$-constellations and to construct their moduli spaces 
$M_\theta$. In a quiver-theoretic context, Kronheimer \cite{Kron89} and 
Sardo-Infirri \cite{Sar-In96a}, \cite{Sar-In96b} have already
considered these moduli spaces and have studied the chamber structure 
in the space $\Pi$ of stability parameters $\theta$, 
where all values of $\theta$ in the same chamber yield the same $M_\theta$. Bridgeland, King and Reid \cite{BKR01} use derived category methods 
to show, in case of arbitrary $G \subseteq \gsl_3(\mathbb{C})$ 
that $G$-$\hilb$ is a crepant resolution of $X$. 
Their method can be used to show that, for any chamber in $\Pi$, $M_\theta$ 
is a crepant resolution, however it yields little information about 
either the structure of the chamber space or the geometry 
of $M_\theta$s. 

  Craw in his thesis conjectured that every projective crepant resolution 
of $X$ can be realised as a moduli space $M_\theta$ of $\theta$-stable
$G$-constellations for some chamber in $\Pi$. A recent paper by  
Craw and Ishii \cite{Craw-Ishii-02} proves this for all abelian 
$G \subset \gsl_3(\mathbb{C})$.

  In this paper, we take a different approach to this issue. Rather
than constructing a resolution as a moduli space of
$G$-constellations, we shall take an arbitrary (not necessarily
projective or crepant) resolution of $X$ and study what families of
$G$-constellations it can parametrise. 

  To start with let $G$ be any finite abelian subgroup of $\gl_n
(\mathbb{C})$ and $Y$ any scheme birational to the quotient space 
$X = \mathbb{C}^n / G$. 
\begin{align*}
\xymatrix{ 
Y \ar[dr]^{\pi} & & \mathbb{C}^n \ar[dl]_{q} \\
& X &
}
\end{align*}

  Let $R$ denote the coordinate ring $\mathbb{C}[x_1, \dots, x_n]$ of
$\mathbb{C}^n$. A \it $(G,\regring)$-module \rm is a
$G$-representation $V$ together with a $G$-equivariant action of
$\regring$.  The categories of finite-length $G$-equivariant coherent
sheaves on $\mathbb{C}^n$ and of $(G,\regring)$-modules are equivalent 
and in this paper we work in the latter category. 

  We would like the families of $G$-constellations which we study
to be related, geometrically, to the space $Y$ which parametrises them. 
That is, we would like to single out a set of `natural' families of 
$G$-constellations on $Y$. For instance, for any point $y \in Y$ we
have its image $\pi(y)$ in $X$ and hence an orbit $q^{-1}(\pi(y))$ of $G$ 
in $\mathbb{C}^n$. On the other hand, a $G$-constellation is a
$G$-equivariant finite-length sheaf and hence is supported on a finite union
of $G$-orbits in $\mathbb{C}^n$. It seems reasonable to ask 
for the $G$-constellation parametrised by $y \in Y$ to be supported, 
set theoretically, precisely on $q^{-1}(\pi(y))$. 

Observe now that, due to dimension considerations, there is only
one $G$-constellation supported at any free orbit of $G$ in $\mathbb{C}^n$, 
up to an isomorphism. This $G$-constellation is precisely 
the structure sheaf $\mathcal{O}_Z$ of $G$-cluster $Z$ given by that orbit. 
Thus $q_* \mathcal{O}_{\mathbb{C}^n}$, 
over any subset $U$ of $X$ such that $G$ acts freely on $q^{-1}(U)$, 
is a unique (up to a twist by a line bundle) family of $G$-constellations
satisfying the wanted property on supports. Observe, that its fiber 
at the generic point of $X$ is the $G$-constellation 
$K(\mathbb{C}^n) \simeq \regrep \otimes K(X)$, which we can think of
as corresponding to the generic orbit of $G$.  As any scheme
birational to $X$ shares its generic point $p_X$, the very least any
natural family should do is to have $p_X$ parametrise a
$G$-constellation isomorphic to $K(\mathbb{C}^n)$. We call such
families \it deformations of the generic orbit of $G$ across $Y$\rm.
We then show (Proposition \ref{prps-tfae}) that this requirement on 
the fiber of the family at the generic point implies much 
stronger naturality properties: for any point $y \in Y$, 
the support of the $G$-constellation it parametrises is indeed
$q^{-1} \pi(y)$, set-theoretically.  Moreover, any such family can 
be $G$ and $R$ equivariantly embedded into the constant sheaf 
$K(\mathbb{C}^n)$ on $Y$.

Now for $G$ abelian, any family of $G$-constellations is a direct sum 
of invertible $G$-eigensheaves. On any scheme $S$, to consider invertible
$\mathcal{O}_S$-submodules of $K(S)$ is to consider Cartier 
divisors on $S$. Therefore in Section \ref{section-valuations} 
we extend the construction of Cartier divisors on $Y$, as 
global sections of $K^*(Y)/\mathcal{O}^*_Y$, by defining a \it $G$-Cartier 
divisor \rm to be a global section of $K^*_G(\mathbb{C}^n)/\mathcal{O}^*_Y$, 
where $K^*_G(\mathbb{C}^n)$ is the group of all non-zero $G$-homogeneous 
rational functions on $\mathbb{C}^n$. 

To make a link with Weil divisors, we make the natural extension of the concept
of the valuation at a prime divisor from $K(Y)$ to $K^*_G(\mathbb{C}^n)$. 
We then define $G$-Weil divisors (Definition \ref{defn-gweil-divisor})
as a subset of $\mathbb{Q}$-Weil divisors on $Y$, in such a way as
to have the correspondence between $G$-Weil and $G$-Cartier
divisors in place when $Y$ is smooth. 

  Now as any deformation $\mathcal{F}$ of the generic orbit embeds into
$K(\mathbb{C}^n)$ as a $(G,\regring)$-submodule, each of its eigensheaves 
$\mathcal{F}_\chi$, together with its embedding into $K(\mathbb{C}^n)$ 
defines a $G$-Cartier divisor and consequently a $G$-Weil divisor
$D_\chi$. Conversely, any set $\{D_\chi\}$, where for each 
$\chi \in G^\vee$ we have one
$\chi$-Weil divisor $D_\chi$, defines an $\mathcal{O}_{Y}$-submodule
$\oplus \mathcal{L}(- D_\chi)$ of $K(\mathbb{C}^n)$. For it to be a
$(G,\regring)$-submodule, and hence a deformation of the generic
orbit, we need the $R$-action on $K(\mathbb{C}^n)$ to restrict down to
it. We show that this 
is precisely equivalent to the condition that for $(f)$ the principal 
divisor of any $G$-homogeneous $f \in \regring$
\begin{align*} D_\chi + (f) - D_{\chi \rho(f)} \geq 0 \end{align*} 
where $\rho(f)$ is the weight of $f$. Now it is clearly sufficient 
for this to be true just for $f = x_1, \dots, x_n$, the basic monomials. 
Thus we establish a 1-to-1 correspondence between deformations of the
generic orbit and sets $\{D_\chi\}_{\chi \in G^\vee}$ of $G$-Weil 
divisors satisfying a finite number of inequalities. 

 It is usual in moduli problems to consider the families up to
equivalence, that is twisting by a line bundle. We show that any 
equivalence class of deformations of the generic orbits contains 
a unique family with $D_{\chi_0} = 0$ in the corresponding divisor set. 
We call such deformations of the generic orbit \it normalized\rm. 
On the other hand, the requirement for the subsheaf $\oplus \mathcal{L}(-D_\chi)$ of
$K(\mathbb{C}^n)$ to be closed under $R$-action can be seen to imply 
that all the eigensheaves $\mathcal{L}(-D_\chi)$ must be, in a certain sense, 
close to each other inside $K(\mathbb{C}^n)$. 
When $D_{\chi_0} = 0$, this allows us to put a precise
bound on how far from $0$, numerically, all the other divisors $D_\chi$ can be. 
Explicitly, we define the set $\{M_\chi\}$ by 
\begin{align*} 
M_\chi = \sum_P (\min_{f \in \regring_\chi} v_P(f)) P
\end{align*}
where $P$ ranges over all prime Weil divisors on $Y$. We show that
$\oplus \mathcal{L}(-M_\chi)$ is a deformation of the generic orbit, 
and in case of $Y$ being $G$-$\hilb$ it is the tautological family 
of $G$-clusters parametrised by $Y$. Then we prove that for 
any normalized deformation of the generic orbit, the corresponding 
divisor set $\{D_\chi\}$ satisfies
\begin{align*} 
M_\chi \geq
D_\chi \geq
- M_{\chi^{-1}}
\end{align*}

In particular, this implies that the number of equivalence classes is
finite as we show that the only non-zero summands of $M_\chi$ are the
exceptional divisors and the proper transforms in $Y$ of images in
$X$ of coordinate hyperplanes of $\mathbb{C}^n$. 

Thus our main result (Theorem \ref{theorem-classification}) is:

\begin{theorem*}[Classification]
Let $G$ be a finite abelian subgroup of $\gl_n(\mathbb{C})$, $X$ be the
quotient of $\mathbb{C}^n$ by the action of $G$ and $Y$ be a resolution
of $X$. Then all deformations of the generic orbit across $Y$, 
up to isomorphism, are of form 
$\oplus_{\chi \in G^\vee} \mathcal{L}(-D_\chi)$, where each 
$D_\chi$ is a $\chi$-Weil divisor and the set $\{D_\chi\}$ satisfies 
the inequalities:
\begin{align*}
D_\chi + (f) - D_{\chi \rho(f)} \geq 0 
\end{align*}
for all $\chi \in G^\vee$ and all $G$-homogeneous $f \in \regring$.
Here $\rho(f)$ is the homogeneous weight of $f$. Conversely for any 
such set $\{D_\chi\}$, $\oplus \mathcal{L}(-D_\chi)$ is a deformation
of the generic orbit.     

Moreover, each equivalence class of families has precisely one family 
with $D_{\chi_0} = 0$. The divisor set $\{D_\chi\}$ corresponding 
to such a family satisfies inequalities 
\begin{align*} M_\chi \geq
D_\chi \geq
- M_{\chi^{-1}}
\end{align*}
where $\{ M_{\chi} \}$ is a fixed divisor set depending only on $G$
and $Y$. In particular, the number of equivalence classes of families 
is finite. 
\end{theorem*}

Throughout the paper we illustrate the proceedings with examples from 
toric geometry, which allows for explicit 
calculations on $Y$ whenever $G$ is abelian. A brief summary of 
the toric setup as applied to our problem is given in Section 
\ref{section-toric-picture}. Then we introduce $Y$ on which all 
of the examples will be calculated: a single toric flop of
$G$-$\hilb$, with $G$ being the cyclic subgroup of $\gl_3(\mathbb{C})$ 
of order $8$ traditionally denoted $\frac{1}{8}(1,2,5)$. 

\bf Acknowledgements: \rm The author would like to express his gratitude to Alastair Craw and Akira
Ishii for many useful discussions on the subject and to Alastair King
for all the insights, corrections and for his tireless support in shaping 
this paper into its present form.

\section{Deformations of the Generic Orbit} \label{section-deform}

\subsection{$G$-Constellations and Families}

Let $G$ be a finite abelian group and let $\givrep$ be an
$n$-dimensional faithful representation of $G$. We identify the symmetric 
algebra $S(\givrep^\vee)$ with the coordinate ring $\regring$ 
of $\mathbb{C}^n$ via a choice of such an isomorphism that 
the induced action of $G$ on $\mathbb{C}^n$ is diagonal. 
By the dual action of $G$ on $\regring$ we shall mean 
the left action given by
\begin{align} \label{eq-gaction}
  g.f(\bv) = f(g^{-1}.\bv)\quad\quad\forall\; \bv \in \mathbb{C}^n
\end{align}

 Corresponding to the inclusion $\regring^G \subset \regring$ of the 
subring of $G$-invariant functions we have the quotient map $q:
\mathbb{C}^n \rightarrow X$, where $X = \spec \regring^G$ is 
the quotient space. This space is generally singular. So we are
typically interested in taking resolutions $\pi: Y \rightarrow X $ of
it.
\begin{align*}
\xymatrix{ 
Y \ar[dr]^{\pi} & & \mathbb{C}^n \ar[dl]_{q} \\
& X &
}
\end{align*}
 The purpose of this paper is to study the way in which $Y$ can
parametrise families of \it $G$-constellations\rm. 

\begin{defn}[\cite{Craw-Ishii-02}] \label{gcon-as-sheaf}
A \tt $G$-constellation \rm is a $G$-equivariant coherent sheaf
$\mathcal{F}$ on
$\mathbb{C}^n$ such that $H^0(\mathcal{F})$ is isomorphic, as a
$\mathbb{C}[G]$-module, to the regular representation $\regrep$. 
\end{defn}

Of course as $\mathcal{F}$ is coherent, it is uniquely determined by
$H^0(\mathcal{F})$ via the $\tilde{\bullet}$ construction (\cite{Harts77}, p.
110). The actions of $G$ and $R$ on $\mathcal{F}$ are
entirely determined by their restrictions to $H^0(\mathcal{F})$. 
In this paper we shall adopt this 
more algebraic point of view, and consider a following class of
objects:

\begin{defn} \label{defn-gcon}
A \tt $(G,\regring)$-module \rm is a $\mathbb{C}[G]$-module $V$
together with an equivariant $\regring$-action, that is 
\begin{align} \label{o-g-equiv}
g.(f.\bv) = (g.f).(g.\bv)
\end{align}
must hold for all $\bf v \rm \in V$, $g \in G$ and all $f \in \regring$. 

A morphism of $(G,\regring)$-modules is a $G$ and $R$ equivariant
linear map of the underlying vector spaces. 

\end{defn}

The functors $\tilde{\bullet}$ and $H^0(\bullet)$ provide an
equivalence between the categories of finite-length coherent $G$-equivariant 
sheaves on $\mathbb{C}^n$ and of $(G,\regring)$-modules, thus
we can can use both concepts interchangeably.

Any $R$-action on $V$ is defined by an element
of $\homm_\mathbb{C}(\regring \otimes_\mathbb{C} V, V)$. As
$\regring = S(\givrep^\vee)$ it is sufficient to consider restrictions 
to $\homm_\mathbb{C}(\givrep^\vee \otimes V, V)$. 
The condition \eqref{o-g-equiv} is precisely equivalent to asking for 
this homomorphism to be $G$-equivariant.  

Conversely, $\alpha \in \homm_G(\givrep^\vee \otimes V, V)$ 
defines an $\regring$-action on $V$ if and only if it satisfies 
\begin{align}\label{eqn-commutator-conditions-gen}
\alpha(v_1 \otimes \alpha(v_2 \otimes v)) = 
\alpha(v_2 \otimes \alpha(v_1 \otimes v))
\end{align}

Thus we see that there exists a one-to-one correspondence between all
the $(G,\regring)$-modules with an underlying $\mathbb{C}[G]$-module
$V$ and the elements of $Z_{\regring, G} \subseteq
\homm_G(\givrep^\vee~\otimes~V, V)$ satisfying the commutator
conditions \eqref{eqn-commutator-conditions-gen}.

Further, it can be seen that the $R$-structures of two isomorphic
$(G,\regring)$-modules on $V$ differ by conjugation by an element of
$\autm_G(V)$. Therefore we have a one-to-one correspondence between
isomorphism classes of $(G,\regring)$-modules with underlying
$\mathbb{C}[G]$-module $V$ and the orbits of $\autm_G(V)$ in
$Z_{\regring, G}$.  

\begin{defn} \label{def-gcon-fam}
A \tt family of $(G,\regring)$-modules parametrised by a scheme $S$ \rm 
is a locally free sheaf $\mathcal{F}$ of $\mathcal{O}_S$-modules 
with $G$ and $\regring$ acting by $\mathcal{O}_S$-linear
endomorphisms, so that 
\begin{align*}
g.(f.s) = (g.f).(g.s)
\end{align*}
for all $g \in G$, $f \in \regring$ and any local section $s$ of
$\mathcal{F}$. 

 We shall say that two families $\mathcal{F}$ and $\mathcal{F}'$ are 
equivalent if there exists an invertible sheaf $\mathcal{L}$ on $S$ 
such that $\mathcal{F}$ is $(G,\regring)$-equivariantly isomorphic 
to $\mathcal{F}' \otimes \mathcal{L}$.
  
 We shall call $\mathcal{F}$ a \tt family of $G$-constellations \rm if
its fiber $\mathcal{F}_{|p}$ at any point $p \in Y$ is a
$G$-constellation.
\end{defn}

  Any sheaf $\mathcal{F}$ with a $G$-action must split into 
$G$-eigensheaves, which are locally-free if $\mathcal{F}$ is. 
In particular, we see that for an abelian $G$ any family 
of $G$-constellations must split as 
$$ \bigoplus_{\chi \in G^\vee} \mathcal{L}_\chi $$
where $G$ acts on each invertible sheaf 
$\mathcal{L}_\chi$ by the character $\chi$. 

  Any free $G$-orbit $Z \subset \mathbb{C}^n$ is a $G$-cluster, 
its structure sheaf $\mathcal{O}_Z$ a $G$-constellation. Considering
$H^0(\mathcal{O}_Z)$ as the fibre of $q_* \mathcal{O}_{\mathbb{C}^n}$
at $x = q(Z) \in X$, we see that over any $U \subset X$ such that $G$ acts
freely on $q^{-1}(U)$, we have a natural family of $G$-constellations
$\mathcal{F} = q_*\mathcal{O}_{\mathbb{C}^n}$. 

  Now consider the generic point $\genpx_$ of $X$. Its pre-image in 
$\mathbb{C}^n$ is the generic point $\genpcn_$, which can be viewed
as the generic orbit of $G$. The fibre of
$\mathcal{O}_{\mathbb{C}^n}$ at $\genpcn_$ is the function field
$K(\mathbb{C}^n)$ and that of $\mathcal{O}_X$ at $\genpx_$ is $K(X) = 
K(\mathbb{C}^n)^G$. The extension $K(\mathbb{C}^n) : K(\mathbb{C}^n)^G$
is Galois, so the Normal Basis Theorem from Galois theory
(\cite{Garl86}, Theorem 19.6) implies that 
$K(\mathbb{C}^n) = \regrep \otimes_{\mathbb{C}} K(X)$.
Thus $K(\mathbb{C}^n)$ is a family of $G$-constellations 
parametrised by a single point-scheme $p_X$. 
Moreover it is natural, in the sense that it is precisely 
the fiber of the natural family $q_*(\mathcal{O}_{\mathbb{C}^n})$ at $p_X$. 
We now proceed to single out a class of families whose fiber at the
generic point is isomorphic to the natural one. 

\begin{defn}\label{defn-deform-gen-orbit}
Let $Y$ be a scheme birational to $X$ and let $p_Y$ denote 
the generic point of $Y$. 
A \tt deformation of the generic orbit of $G$ across $Y$ \rm 
is a family of $G$-constellations parametrised by $Y$ equipped with
a $(G,\regring)$-equivariant isomorphism $$ \iota :\;  \mathcal{F}_{|p_Y} 
\xrightarrow{\sim} K(\mathbb{C}^n) =  (\pi^* q_*
\mathcal{O}_{\mathbb{C}^n})_{|p_Y} $$ 
\end{defn}

 We now show that, in fact, any family which agrees with the natural
one at the generic point must agree with it wherever $G$ acts freely.  

\begin{prps} \label{prps-tfae}
Let $\pi: Y \rightarrow X$ be a birational morphism and let $F$ be a family 
of $G$-constellations on $Y$. Then the following are equivalent:
\begin{enumerate}
\item \label{item-genpt} There exists an isomorphism 
\begin{align} \label{eqn-deform-genorb}
\mathcal{F}_{|\genpy_} \simeq K(\mathbb{C}^n)
\end{align}
which makes $\mathcal{F}$ into a deformation of the generic orbit of
$G$ across $Y$.
\item \label{item-embedding} There exists a $(G,\regring)$-equivariant
embedding $\iota': \mathcal{F} \hookrightarrow K(\mathbb{C}^n)$, where
$K(\mathbb{C}^n)$ is considered as a constant sheaf of $(G,\regring)$-modules 
on $Y$. 
\item \label{item-actions}
For any open $U \subseteq Y$, $s \in \mathcal{F}(U)$ and $f \in R^G$ we have
\begin{align} \label{eqn-r-ou-action}
f.s = f s 
\end{align}
where on the left-hand side $f$ acts as an element of $R$ and on the
right-hand side as a section of $\mathcal{O}_Y$, via the inclusion
$\mathcal{O}_X \hookrightarrow \pi_* \mathcal{O}_Y$.  
\item \label{item-extension} For any open $U \subset X$ such that 
$G$ acts freely on $q^{-1}U$, 
\begin{align} \label{eqn-extension}
\mathcal{F}|_{\pi^{-1}U}  \simeq 
\pi^* q_*\mathcal{O}_{\mathbb{C}^n}|_{\pi^{-1}U} \otimes \mathcal{L}
\end{align}
for some invertible sheaf $\mathcal{L}$ on $\pi^{-1}U$. 
\end{enumerate}
\end{prps}

Before tackling this proposition, we prove a useful lemma, which provides 
a nice geometrical interpretation of the condition \eqref{eqn-r-ou-action}. 

\begin{lemma} \label{lemma-fibers-iso}
Let $\mathcal{F}$ be a family of $G$-constellations on $Y$ satisfying
$\eqref{eqn-r-ou-action}$. Then for any $p \in Y$ we have a 
scheme-theoretic inclusion 
\begin{align} \label{eqn-support-inclusion}
\supp \mathcal{F}_{|p} \subseteq q^{-1}\pi(p) 
\end{align}
Moreover, set-theoretically we have equality. Further, if $G$ acts
freely on $q^{-1}(p)$, we have
$$ \mathcal{F}_{|p} \simeq (\pi^* q_* \mathcal{O}_{\mathbb{C}^n})_{|p} $$
as $G$-constellations. 
\end{lemma}
\begin{proof}
Given an arbitrary $G$-constellation $V$, the support of $V$ is the
vanishing set of the ideal $\ann_\regring V \subset \regring$. On 
the other hand $q^{-1}\pi(p)$ is the vanishing of the ideal in 
$R$ generated by $\mathfrak{m}_{\pi p} \in \regring^G$. 
So scheme-theoretically \eqref{eqn-support-inclusion} is equivalent to 
$$ \ann_{R^G} k(\pi p) \subset \ann_{R} \mathcal{F}
\otimes_{\mathcal{O}_Y} k(p) $$
which follows immediately from \eqref{eqn-r-ou-action}. 

To show the set-theoretic equality, we observe from \eqref{o-g-equiv} that 
the ideal $\ann_R \mathcal{F}_p$ is $G$-invariant, and
so, set-theoretically $\supp \mathcal{F}_{|p}$ is a union of
$G$-orbits in $\mathbb{C}^n$. But \eqref{eqn-support-inclusion} now
implies that it is contained in a single orbit: the closed points of 
$q^{-1}\pi(p)$. Therefore we have equality. 

For the last bit, we observe that $\mathcal{F}_{|p}$ is a finite
length sheaf on $\mathbb{C}^n$ and so splits as a direct sum 
$$ \bigoplus_{x \in \supp \mathcal{F}_{|p}} (\mathcal{F}_{|p})_{|x} $$
of its fibers at each closed point in its support. But as $G$ acts
freely on $q^{-1}\pi(p)$, the size of the orbit is $|G|$. Since
this is also the dimension of $\mathcal{F}_{|p}$, each 
$(\mathcal{F}_{|p})_{|x}$ must be 1-dimensional and hence 
$$ \mathcal{F}_{|p} = \bigoplus_{x \in q^{-1}\pi(p)}
(\mathcal{O}_{\mathbb{C}^n})_{|x} \simeq (\pi^* q_*
\mathcal{O}_{\mathbb{C}^n})_{|p} $$
\end{proof}

\begin{proof}[Proof of Proposition \ref{prps-tfae}]

$\ref{item-extension} \Rightarrow \ref{item-genpt}$ is obtained by
considering the restriction of the isomorphism $\eqref{eqn-extension}$ 
to stalks at $\genpy_$. 

$\ref{item-genpt} \Leftrightarrow \ref{item-embedding}$: consider
the sheaf $\mathcal{F} \otimes_{\mathcal{O}_Y} K(Y)$. 
On any open $U$ where $\mathcal{F}$ is a free
$\mathcal{O}_Y$-module, $ \mathcal{F} \otimes_{\mathcal{O}_Y} K(Y) $
is the constant sheaf $\mathcal{F}_{p_Y}$ for which we have the
$(G,\regring)$-equivariant isomorphism \eqref{eqn-deform-genorb} 
to the constant sheaf $K(\mathbb{C}^n)$. A sheaf constant on an open cover 
must be constant as $Y$ is irreducible.  Now the natural map 
$\mathcal{F} \hookrightarrow \mathcal{F} \otimes K(Y)$ 
becomes the requisite embedding.

$\ref{item-embedding} \Rightarrow \ref{item-actions}$ is immediate
because $K(\mathbb{C}^n)$, as a $(G,\regring)$-module clearly satisfies 
\eqref{eqn-r-ou-action}.

So we are left with proving 
$\ref{item-actions} \Rightarrow \ref{item-extension}$. 

We begin with a local version: if $p \in \pi^{-1}(U) \subset
Y$, then $\mathcal{F}_p \simeq (\pi^* q_*
\mathcal{O}_{\mathbb{C}^n})_p$, that is the stalks at $p$ are
$(G,\regring)$-equivariantly isomorphic. 

Now $(\pi^* q_* \mathcal{O}_{\mathbb{C}^n})_p$ (which we can write as 
$\regring \otimes_{\regring^G} \mathcal{O}_{Y,p}$) is a free 
$\mathcal{O}_{Y,p}$-module of rank $|G|$. This is because $G$ acting 
freely on $q^{-1} \pi(p)$ implies that the quotient map $q$ is flat 
and $|G|$-to-one at $\pi(p)$. $\mathcal{F}_p$ is also a free 
$\mathcal{O}_{Y,p}$-module of rank $|G|$, because $\mathcal{F}$ is 
a family of $G$-constellations. Therefore we can consider the 
determinant of any $(G,\regring)$-equivariant $\mathcal{O}_{Y,p}$-morphism 
between the two, and it would suffice to find a morphism whose
determinant is invertible. 

Consider the map $\theta:\; (\pi^* q_* \mathcal{O}_{\mathbb{C}^n})_p
\rightarrow \mathcal{F}_p$ defined by
\begin{align}
m \otimes f \rightarrow m.(f s_0) \quad m \in \regring,\;f \in \mathcal{O}_{Y,p}
\end{align}
where $s_0$ is a fixed choice of any $\mathcal{O}_{Y,p}$-generator of 
the $\chi_0$-eigenspace of $\mathcal{F}_p$. 

This map is a well-defined $\mathcal{O}_{Y,p}$-module map, that is it
descends from the set-theoretic product $\regring \times \mathcal{O}_{Y,p}$  
to the tensor product, precisely because both $\mathcal{F}_p$ and 
$\regring \otimes \mathcal{O}_{Y,p}$ satisfy \eqref{eqn-r-ou-action}.
It is $G$-equivariant because $1 \mapsto s_0$ ensures that
$\chi_0$-eigenspace maps to $\chi_0$-eigenspace and \eqref{o-g-equiv}
forces the rest. Finally not only $\theta$ is defined to be
$\regring$-action equivariant, but the reader can verify that it is the
unique element of $\homm_{(G,\regring)}(\regring \otimes
\mathcal{O}_{Y,p}, \mathcal{F}_p)$ which maps $1$ to $s_0$.
Note that in particular, this shows that 
\begin{align*} 
\homm_{(G,\regring)}(\regring \otimes \mathcal{O}_{Y,p},
\mathcal{F}_p) \simeq (\mathcal{F}_{p})_{\trch_} \simeq
\mathcal{O}_{Y,p} \quad\quad (\dagger)
\end{align*}
$\theta$ is a $(G,\regring)$-equivariant morphism. It descends to
the $(G,\regring)$-equivariant morphism  
$$\overline{\theta}: (\pi^* q_* \mathcal{O}_{\mathbb{C}^n})_{|p}
\rightarrow \mathcal{F}_{|p}$$ on fibers. Similarly to $(\dagger)$, 
\begin{align*} 
\homm_{(G,\regring)}((\pi^* q_* \mathcal{O}_{\mathbb{C}^n})_{|p}, 
\mathcal{F}_p) \simeq \mathbb{C}
\end{align*}
i.e. all $(G,\regring)$-equivariant morphisms between the two are 
scalar multiples of each other. Since by Lemma \ref{lemma-fibers-iso}, 
the two fibers are $(G,\regring)$-equivariantly isomorphic, we have
that unless $\overline{\theta}$ is a zero map, it is an isomorphism. 
But it maps $[1]$ to $[s_0]$, and the latter can not be $0$ 
by the choice of $s_0$. So $\det \overline{\theta} \neq 0$ 
impying that $\det \theta \in \mathcal{O}^*_{Y,p}$, as required. 

The isomorphisms on stalks give isomorphisms 
$\theta_i : \regring \otimes_{\regring^G} \mathcal{O}_{U_i} \rightarrow
\mathcal{F}|_{U_i}$ on an open 
cover $\{ U_i \}$ of $U$, as both sheaves are locally free and of finite rank. 
Then on each intersection $U_i \cap U_j$, $\theta_i \circ \theta_j^{-1}$
is a $(G,\regring)$-automorphism of $\regring \otimes_{\regring^G}
\mathcal{O}_{U_i \cap U_j}$. Any such, by an argument identical to 
$(\dagger)$, is a multiplication by an element of 
$\mathcal{O}^*_{U_i \cap U_j}$,which concludes the proof.
\end{proof}

 For the rest of this paper, we shall concern ourselves only with those 
families of $G$-constellations which are deformations of the generic
orbit.  

 Observe that the map $\iota:\;  \mathcal{F}_{|p_Y} \xrightarrow{\sim} 
K(\mathbb{C}^n) =  (\pi^* q_* \mathcal{O}_{\mathbb{C}^n})_{|p_Y} $
uniquely determines the embedding $\iota':\; \mathcal{F}
\hookrightarrow K(\mathbb{C}^n)$. The notion of the isomorphism of 
deformations demands for the $(G,\regring)$-equivariant sheaf isomorphism
$\theta:~\mathcal{F}~\rightarrow~\mathcal{F}'$ to have its restriction 
to stalks at $p_Y$ form a commutative triangle with maps
$\iota_{\mathcal{F}}$ and $\iota_{\mathcal{F}'}$ for $\mathcal{F}$ and 
$\mathcal{F}'$ to be isomorphic as deformations of the generic orbit. 
Consequently $\theta$ itself must form a commutative triangle with 
$\iota'_{\mathcal{F}}$ and $\iota'_{\mathcal{F}'}$, in particular 
images of $\mathcal{F}$ and $\mathcal{F}'$ in $K(\mathbb{C}^n)$ 
must coincide. Thus isomorphism classes of deformations of 
the generic orbit are precisely in one-to-one correspondence 
with deformations of the generic orbit which are subsheaves 
of $K(\mathbb{C}^n)$.

\section{Line bundles and $G$-Cartier divisors} \label{section-valuations}

 As we deal with families of $G$-constellations which are subsheaves
of $K(\mathbb{C}^n)$, it would be useful to have a language similar to 
that of the Cartier divisors to describe the invertible 
sub-$\mathcal{O}_Y$-modules of $K(\mathbb{C}^n)$ with non-trivial 
$G$-action. 
 In this section we shall extend the familiar construction of Cartier
divisors using the larger group of non-zero $G$-homogeneous rational
functions, which we shall denote by $\kgc_$, instead of the group of 
non-zero invariant rational functions $K^*(Y)$. 

\begin{defn}
 We shall say that a rational function $f \in K(\mathbb{C}^n)$ is \tt
$G$-homogeneous of weight $\chi \in G^\vee$ \rm if 
such that 
\begin{align} \label{eq-ghom-fn}
g.f = \chi(g^{-1}) f \quad \forall\; g \in G 
\end{align}

 We shall denote by $K_\chi(\mathbb{C}^n)$ the subset of $K(\mathbb{C}^n)$ of
$G$-homogeneous elements of a specific weight $\chi$ 
and by the $K_G(\mathbb{C}^n)$ the subset of $K(\mathbb{C}^n))$ of
all the $G$-homogeneous elements. 
We shall use $\regring_\chi$ and $\regring_G$ to mean
$\regring \cap K_\chi(\mathbb{C}^n)$ and $\regring \cap
K_G(\mathbb{C}^n)$ respectively. 
\end{defn} 

 The choice of a sign in this definition is motivated 
as follows: we want a function $p \in \regring$ to be $G$-homogeneous
of weight $\chi \in G^\vee$ if $p(g.v) = \chi(g) p(v)$ for any $g \in G$ and 
$v \in \mathbb{C}^n$. E.g. usual concept of a homogeneous polynomial,
whose degree, an integer number, is precisely its weight as a character 
of $\mathbb{C}^*$ acting diagonally on $\mathbb{C}^n$. 
In view of \eqref{eq-gaction}, this means we have to have
$\chi(g^{-1})$ instead of $\chi(g)$ in \eqref{eq-ghom-fn}.  

 Now consider $\kgc_$, the invertible elements of
$K_G(\mathbb{C}^n)$. Using the fact that $K(Y) = K(X) = K(\mathbb{C}^n)^G$, 
we have a short exact sequence of multiplicative groups: 

\begin{align}
 1 \rightarrow K^*(Y) \rightarrow \kgc_ \rightarrow G^\vee \rightarrow 1
\end{align}
What makes this enlargement of $K^*(Y)$ useful is that we can still
define a valuation of a $G$-homogeneous rational function at a prime
Weil divisor.

\begin{defn}
 Let $D \subset Y$ be a prime Weil divisor on $Y$. Given any 
$f \in \kgc_$, we choose any $n \in \mathbb{Z}$ such that $f^n$ is invariant, 
i.e. $f^n \in K(Y)$. For instance, $n = |G|$. Then we define
\begin{align}
v_D(f) = \frac{1}{n}v_D(f^n) \in \mathbb{Q}
\end{align}
 where $v_D(f^n)$ is the ordinary valuation of $f^n$ in the local ring
 $\mathcal{O}_{D,Y}$ of the generic point of $D$. This is well-defined
since for any $g \in K(Y)$, we have $v_D(g^k) = kv_D(g)$.
\end{defn}

 In what follows, we shall write $$ \{n\} = n - [ n ] $$
for the fractional part of $n \in \mathbb{Q}$. Generally, 
the valuations defined above are $\mathbb{Q}$-valued. However, 
if $f$ and $g$ in $\kgc_$ are both $\chi$-homogeneous, 
then $f/g$ is $G$-invariant and hence for any Weil divisor $D$ on $Y$, 
$v_D(f) - v_D(g) \in \mathbb{Z}$. Therefore the fractional part of
${ v_D(f) }$ is independent of the choice of $f$ in $\kgchi_$.  

\begin{defn} \label{defn-fract-part}
We define $v(D,\chi)$ to be the number $\{v_D(f)\} \in \mathbb{Q}$, where
$f$ is any element of $\kgchi_$. 
\end{defn}

  We can now replicate, almost word-for-word, the definitions in
\cite{Harts77}, pp. 140-141. 

\begin{defn}
  A \tt $G$-Cartier divisor \rm on $Y$ is a global section of the sheaf 
of multiplicative groups $\kgc_ / \mathcal{O}^*_Y$, i.e. the quotient
of the constant sheaf $\kgc_$ on $Y$ by the sheaf $\mathcal{O}^*_Y$ of
invertible regular functions. 

  As usual, such a section can be described by a choice of an open cover 
$\{U_i\}$ of $Y$ and functions $\{f_i\} \subseteq \kgc_$ such 
that $f_i / f_j \in \Gamma(U_i \cap U_j, \mathcal{O}^*_Y)$. Observe
that, as their ratios are invariant, the $f_i$ must all be homogenous 
of the same weight $\chi \in G^\vee$.  In such a case, we further say that
the divisor is \tt $\chi$-Cartier\rm. 
\end{defn}

  As with ordinary Cartier divisors, a $G$-Cartier divisor is said to be
principal if it lies in the image of the natural map 
$\kgc_ \rightarrow \kgc_ / \mathcal{O}^*_Y$ and two divisors are
said to be linearly equivalent if their difference is principal. 

  However when defining a corresponding enlargement of 
the group of Weil divisors, we have to be a little bit careful. 

\begin{defn} \label{defn-gweil-divisor}
A \tt $\chi$-Weil divisor \rm on Y is a finite sum $\sum q_i D_i$ (where $q_i
\in \mathbb{Q}$) of prime Weil divisors on $Y$, such that 
\begin{align} \label{eqn-gweil-cond}
q_i - v(D_i, \chi) \in \mathbb{Z}
\end{align}
for all $i$. 

  We shall further use the term \tt $G$-Weil divisor \rm to refer to
all $\chi$-divisors for any $\chi \in G^\vee$. 
\end{defn}

\begin{defn}
For any $f \in \kgc_$, we define \tt the principal $G$-Weil 
divisor of $f$ \rm to be 
\begin{align*}
(f) = \sum v_P(f) P
\end{align*}
with the sum taken over all prime Weil divisors $P$ on $Y$. This
sum is finite as $f^{|G|}$ is a regular function on $Y$ and hence has 
non-zero valuations only on finitely many prime divisors.
\end{defn}
 
 Given any $\chi, \chi' \in G^\vee$, we can see that, for 
any prime divisor $D$,
$$ v(D,\chi) + v(D,\chi') - v(D, \chi \chi') \in \mathbb{Z} $$
as it is equal to the valuation at $D$ of an invariant function. 
Hence $G$-Weil divisors form an additive group. 
We define two $G$-Weil divisors to be linearly equivalent if their difference
is principal and a divisor $\sum q_i D_i$ to be effective if
all $q_i \geq 0$. 

Recall (\cite{Harts77}, Proposition 6.11) that there is an injective 
homomorphism from the group of Cartier divisors to the group of Weil
divisors which is an isomorphism when $Y$ is smooth. The definition 
extends naturally to an injective homorphism from the group 
of $G$-Cartier divisors to the group of $G$-Weil divisors, but some care 
needs to be taken to show that it is surjective when $Y$ is smooth.

\begin{defn}
Define the map $\phi$ from the group of $G$-Cartier divisors to the group of 
$G$-Weil divisors on $Y$ by 
$$\{(f_i, U_i)\} \mapsto \sum k_D D$$ where 
the sum is taken over all prime Weil divisors $D$ on $Y$ and 
$k_D = v_D(f_i)$ for any $f_i$ such that $U_i \cap D$ is not empty.
Once again the sum is finite, as each $f_i$ has non-zero valuation
only on finitely many prime Weil divisors.
\end{defn}

\begin{prps} \label{prps-cartier-weil}
 Let $\phi$ be the injective homomorphism defined above. If $Y$ is smooth, 
then $\phi$ is an isomorphism.
\end{prps}

\begin{proof}
  We need surjectivity. So suppose we have a $\chi$-Weil divisor $D$ on $Y$. 
Take any $f \in \kgchi_$. Then $D - (f)$ is an ordinary Weil 
divisor and as $Y$ is smooth, it has a Cartier divisor $\{(U_i, g_i)\}$ 
corresponding to it as before. Then $\{(U_i, g_i f)\}$ is 
the $\chi$-Cartier divisor which $\phi$ maps to $D$. 
\end{proof}

  The point of introducing $G$-Cartier divisors is that they correspond 
to invertible sheaves which carry a $G$-action in the same way that
ordinary Cartier divisors correspond to the ordinary invertible sheaves.

  Indeed consider $D$, the $\chi$-Cartier divisor on $Y$ specified 
by a collection $\{(U_i, f_i)\}$ where $U_i$ form 
an open cover of $Y$ and $f_i \in
\kgchi_$. We define an invertible sheaf $\mathcal{L}(D)$ on $Y$ as
the sub-$\mathcal{O}_Y$-module of $K(\mathbb{C}^n)$ generated by
$f_i^{-1}$ on $U_i$. Observe that we have an action of 
$G$ on $\mathcal{L}(D)$, restricted from the one on 
$K(\mathbb{C}^n)$, and it acts on every section by 
the character $\chi$.  

\begin{prps}
  The map $D \rightarrow \mathcal{L}(D)$ gives an isomorphism between 
the group $\gcl_$ of $G$-Cartier divisors up to linear equivalence and
the group $\gpic_$ of invertible $G$-sheaves on $Y$.
\end{prps}
\begin{proof}
A standard argument from \cite{Harts77}, Corollary 6.15, shows that it is 
an injective homomorphism. To show that it is an isomorphism, we need 
to be able to embed any invertible $G$-sheaf $\mathcal{L}$, with 
$G$ acting by some $\chi \in G^\vee$, as a sub-$\mathcal{O}_Y$-module 
into $K(\mathbb{C}^n)$. 
  
  Given such $\mathcal{L}$, we consider the sheaf 
$\mathcal{L} \otimes_{\mathcal{O}_Y} K(Y)$. On every open set $U_i$ where 
$\mathcal{L}$ is trivial, it is $G$-equivariantly isomorphic to the constant 
sheaf $K_\chi(\mathbb{C}^n)$. On an irreducible scheme a sheaf
constant on an open cover is constant itself, so as $Y$ is irreducible we have 
$\mathcal{L} \otimes_{\mathcal{O}_Y} K(Y) \simeq K_\chi(\mathbb{C}^n)$ and a 
particular choice of this isomorphism gives the necessary embedding as
\begin{align*}
\mathcal{L} \rightarrow \mathcal{L} \otimes_{\mathcal{O}_Y} K(Y)
\simeq K_\chi(\mathbb{C}^n) \subset K(\mathbb{C}^n)
\end{align*}
\end{proof}

A curious thing about $G$-divisors and valuations of $G$-homogeneous
functions is the fact that on the quotient space $X$ every prime 
Weil divisor is a principal divisor of some $G$-homogeneous function.
In particular, every $G$-Weil divisor is $G$-Cartier.  

\begin{prps} \label{prps-weil-principal}
 Let $P$ be a prime Weil divisor on $X$. Then there exists an
$f \in \regring^*_{G}$ such that $P = (f)$, that is
\begin{align*}
v_D(f) = 
\begin{cases}
1, \text{ when } D = P \\
0, \text{ when } D \neq P
\end{cases}
\end{align*}
for any prime divisor $D$ on $Y$.
\end{prps}

\begin{proof}
  Let $I_P \subset \regring^G$ be the prime ideal of height 1 corresponding to
$P$. Consider the ring extension $\regring^G \subseteq \regring$. By a classical
result of Emmy Noether (\cite{Bens94}, Theorem 1.3.1), this extension 
is integral. This then implies (\cite{Mats86}, Theorem 9.3) that there 
exists a prime ideal $I'$ of height 1 in $\regring$ lying over $I_P$, that 
is $I_P = I' \cap \regring^G$ and that every other prime ideal lying over 
$I_P$ is conjugate to $I'$ by an element of $G$. As $\regring$ is an
UFD, every prime ideal of height one is principal and so there exists 
some $y' \in \regring$ such that $I_P = (y') \cap \regring^G$.  

  So take $g_0 = 1, g_1, \dots, g_k \in G$ to be such that the principal ideals 
$(y'), (g_1.y'), \dots, (g_k.y')$ are all the distinct prime ideals lying 
over $I_P$. Then we claim that $y = \prod g_i.y'$ is a $G$-homogeneous 
function and that $I_P = (y) \cap R^G$. Indeed, $(h.y) = \cap
((hg_i).y')$. The ideals $((hg_i).y')$ are all distinct prime ideals
lying over $I_P$ and therefore
$$ (h.y) = \cap ((hg_i).y') = \cap (g_i.y') = (y) $$ 
which implies $h.y \in \mathbb{C}^* y$.  For the second claim, observe that 
$I_P = g_i.I_P = (g_i.y') \cap R^G$ for all $i$. Consequently 
$I_P = (\cap (g_i.y') ) \cap \regring^G = (y) \cap \regring^G$.  
 
  Thus we have $I_P = (y) \cap R^G$. Note that $(y)$ is precisely the
vanishing ideal of the pre-image of $P$ in $\mathbb{C}^n$. Now let $k$ 
be the ramification index of the valuation ring extension 
$R^G_{I_P} \subset \regring_{(y)}$. Then for any $w \in
K(\mathbb{C}^n)^G$ we have $v_P(w) = \frac{1}{k} v_{(y)}(w)$, which 
immediately extends to the $\mathbb{Q}$-valued  valuation
$v_P(w)$ of any $G$-homogeneous $w \in \kgc_$. In particular,  we see
that $v_P(y) = \frac{1}{k}$. Now take any other prime divisor $D$ on
$Y$. We have $I_D = (u) \cap R^G$ for some prime $u \in \regring$. If
now $v_{D}(y) \neq 0$, then as $y$ is regular we have $y \in (u)$ and
so $g_i.y' \in (u)$ for some $i$. Then $(u) = (g_i.y)$ and $D = P$.  
  
   Now taking $f = y^k$ finishes the proof. 

\end{proof}

 In the course of the proof of Proposition \ref{prps-weil-principal}, 
we see that the valuations of $G$-homogeneous functions are actually 
non-integer only at ramification divisors of $q$. We now contemplate 
along which actual divisors the ramification can occur.    

\begin{prps} \label{prps-ramification-hyperplanes}
  There are only finitely many prime divisors $P$ on 
$X$ with ramification index greater than $1$. 
More precisely, if we write the ideal of each such $P$ as
$(y) \cap R^G$ for $y \in \regring^*_{G}$ as per Proposition 
\ref{prps-weil-principal}, then we will have at most one $y$ of weight 
$\chi$ for each character $\chi \in G^\vee$.  

  Explicitly, the ramification can only occur along the images of coordinate 
hyperplanes $(x_1)$, \dots, $(x_n)$ of $\mathbb{C}^n$ and in the case of 
$G \subset \gsl_n(\mathbb{C})$ ramification never occurs at all. 

\end{prps}

\begin{proof}

  For each character $\chi \in G^\vee$ fix a $G$-homogeneous function 
$f_\chi \in R$ of weight $\chi$. We further demand that it
is minimal such, in a sense that no element of $R^G$ other than $1$ divides it.
We shall now show that ramification could only occur along one of the 
$(f_\chi) \cap R^G$ and only when $f_\chi$ is the unique function
satisfying these conditions. 

  To see it, take any prime divisor $P$ on $X$. 
Write $I_P =  (y) \cap R^G$ for $y \in \regring^*_{G}$ as per Proposition 
\ref{prps-weil-principal}. Unless
$f_\chi \in (y)$, $v_{(y)}(f_\chi) = 0$ and hence 
$v_{(y)}(\frac{y}{f_\chi}) = 1$ and so there is no ramification along
$P$. But if $f_\chi \in (y)$ then minimality condition forces $f_\chi = y$. 
 
  Explicitly, when $G$ is abelian we know that the character map 
$\rho:~\mathbb{Z}^n~\rightarrow~G^\vee$ is surjective 
(see Section \ref{section-toric-basics}, \eqref{seq-monom-char}). 
Given a character $\chi \in G$, there exists $m \in \mathbb{Z}^n$
such that $x^m = \prod x_i^{m_i}$ is $G$-homogeneous of
weight $\chi$. Then above implies that 
ramification can only occur along $(y) \cap R^G$ if $y$ is monomial. 
But recalling proof of Proposition \ref{prps-weil-principal}, 
$y = \prod g_i.y'$ where $y'$ is prime. This implies $y'$ must be one 
of the basic monomials $x_i$. 

  In case when $G \subseteq \gsl_n(\mathbb{C})$, we know that 
$x_1 \dots x_n$ is invariant. As $v_{(x_i)}(x_1 \dots x_n) = 1$, 
there is no ramification along any of $(x_i) \cap R^G$ either. 

\end{proof}

 Propositions \ref{prps-weil-principal} and \ref{prps-ramification-hyperplanes}
have an immediate corollary in terms of the numbers $v(P,\chi)$ on $X$.  

\begin{cor} \label{prps-fract-parts-zero}
For any $P$, a prime Weil divisor on $X$
which is not a ramification divisor of $q$,
and $\chi \in G^\vee$, there exists a monomial $m \in \regring_\chi$
such that $v_P(m) = 0$. Consequently
\begin{align*}
v(P,\chi) = 0
\end{align*}
\end{cor}

\begin{proof}
 Unless $P = (x_i) \cap \regring^G$, one can take $m$ to be any monomial 
in $\regring$ of weight $\chi$. 
 If $P = (x_i) \cap \regring^G$, then, unless there is ramification at $P$, 
there exists a $p \in \regring^G$ whose valuation at $(x_i)$ in $\mathbb{C}^n$ 
is $1$. Note that we can take $p$ to be monomial by considering its
monomial summands. Then $\frac{p}{x_i} \in \regring_{\chi^{-1}}$ and 
$v_P(\frac{p}{x_i}) = 0$, so we can take $m = \frac{p}{x_i}^{|G|-1}$. 
\end{proof}

  Let us look at some concrete examples of the ramification occuring 
and not occurring. 

\begin{exmpl}

  First consider $G = \frac{1}{3}(1,2)$, 
the group of $3$rd roots of unity embedded into $\gsl_2(\mathbb{C})$ by 
$$\xi \mapsto \left(
\begin{matrix}
\xi^{1} & \\
& \xi^{2}
\end{matrix} \right)
$$ If we write $\chi_k$ for the character 
of $G$ given by $\xi \mapsto \xi^k$, then $x$ is of weight $\chi_1$
and $y$ of weight $\chi_2$. 

  Let $P$ be the image in $X$ of the hyperplane $x =0$. It is a prime Weil 
divisor (but not a Cartier one) given by $(x^3, xy) = (x) \cap R^G$. 
$v_{(x)}(xy) = 1$, so there is no ramification. And consequently, 
$v_P(x) = v_{(x)}(x) = 1$ as $x^3 = (xy)^3 y^{-3}$.

  Now take $G = \frac{1}{4}(1,2)$. Then the divisor $P$ is 
given by $(x^4, x^2 y)$. So we see that 
index of ramification is $v_{(x)}(x^2 y) = 2$ and correspondingly 
$v_P(x) = \frac{1}{2} v_{(x)}(x) = \frac{1}{2}$. 
\end{exmpl}

\begin{cor} \label{cor-except-nonzero}
Let $\pi: Y \rightarrow X$ be a resolution and $P$ a prime Weil 
divisor on $Y$, which is neither exceptional nor a proper
transform of a ramification divisor of $q$ in $X$. 
Then for any $\chi \in G^\vee$ there exists 
$m \in R_\chi$ such that $v_P(m) = 0$, implying 
\begin{align*}
v(P,\chi) = 0
\end{align*}
\end{cor}

\begin{proof}

  This is a straightforward consequence of Corollary
\ref{prps-fract-parts-zero}. Consider $P' = \pi(P)$, the image of 
$P$ in $X$. Unless $P$ is exceptional, $P'$ is a prime Weil divisor
on $X$. Its generic point lies in the open set on which the
resolution map is an isomorphism, which implies that for 
any $f \in K(\mathbb{C}^n)$, $v_P(f) = v_{P'}(f)$. Now Corollary  
\ref{prps-fract-parts-zero} gives the result.

\end{proof}

\section{Toric Picture} \label{section-toric-picture}

\subsection{Basics} \label{section-toric-basics}
 In this section we give a brief exposition of the necessary toric
background and then translate some of the results of Section
\ref{section-valuations} into the toric language. 

A more thorough exposition of toric geometry in general can be
found in \cite{Dani78} and of toric geometry as related to 
quotient singularities in \cite{ItReid96}.

Consider the maximal torus $(\mathbb{C}^*)^n \subset
\gl_n(\mathbb{C})$ containing $G$. We have an exact 
sequence of abelian groups:

\begin{align}\label{seq-g-torus}
\xymatrix{
0 \ar[r] &
G \ar[r] &
(\mathbb{C}^*)^n \ar[r] &
T \ar[r] &
0
}
\end{align}
where $T$ is the quotient torus which acts on the quotient space $X$.

By applying $\homm(\bullet, \mathbb{C}^*)$ to $\eqref{seq-g-torus}$ we obtain
an exact sequence

\begin{align} \label{seq-monom-char}
\xymatrix{
0 \ar[r] &
M \ar[r] &
\mathbb{Z}^n \ar[r]^{\rho} &
G^\vee \ar[r] &
0
}
\end{align}
where $\mathbb{Z}^n$ is thought of as the lattice of exponents of Laurent
monomials. Thus given $m = (k_1, \dots, k_n) \in \mathbb{Z}^n$ we
shall write $x^{m}$ for $x_1^{k_1} \dots x_n^{k_n}$. $M$ is the
sublattice in $\mathbb{Z}^n$ of (exponents of) $G$-invariant Laurent
monomials. 

Note that each Laurent monomial is a $G$-homogeneous function and 
$\rho$ is precisely the weight map, that is 
$x^{m}(g.\bv) = \rho(m)(g)\; x^{m}(\bv)$ for any $\bv \in \mathbb{C}^n$. 

Applying $\homm(\bullet, \mathbb{Z})$ to $\eqref{seq-monom-char}$ we
obtain 
\begin{align*}
\xymatrix{
0 \ar[r] &
\lbar \ar[r] &
L \ar[r]&
\ext^1(G^\vee, \mathbb{Z}) \ar[r] &
0
}
\end{align*}
where we write $\lbar$ for the dual latice of $\mathbb{Z}^n$, $L$
for the dual of $M$ and note that $\homm(G^\vee, \mathbb{Z}) =
0$ as $G^\vee$ is finite and 
$\ext^1(\mathbb{Z}^n, \mathbb{Z}) = 0$ as $\mathbb{Z}^n$ is free. 

Thus we see that $L/\lbar \simeq \ext^1(G^\vee, \mathbb{Z})$. 
Taking an injective resolution of $\mathbb{Z}$
$$ 0 \rightarrow \mathbb{Z} \rightarrow \mathbb{Q} \rightarrow
\mathbb{Q}/\mathbb{Z} \rightarrow 0$$ we see that $\ext^1(G^\vee, \mathbb{Z})
\simeq \homm(G^\vee, \mathbb{Q}/\mathbb{Z})$ as $\homm(G^\vee,
\mathbb{Q})
= 0$. Now a choice of a map $\mathbb{Q}/\mathbb{Z} \rightarrow
\mathbb{C}^*$ which is equivalent to a simultaneous choice of a
primitive $n$-th root of unity for all $n \in \mathbb{N}$, would
give us $$L/\lbar \simeq \homm(G^\vee, \mathbb{C}^*) = G$$
allowing us to identify points in $L/\lbar$ with 
the elements of the group. 

Tautologically, we have a $\mathbb{Z}$-valued pairing between 
$M$ and $L$. This pairing extends naturally to a $\mathbb{Q}$-valued
pairing between $\mathbb{Z}^n$ and $L$. For the purposes of the
exposition to follow, it will be convenient to think of elements of
$L$ as functions on the monomial lattices $M \hookrightarrow
\mathbb{Z}^n$.  Henceforth, given $l \in L$ and $m \in \mathbb{Z}^n$,
we shall write $l(m)$ to denote the pairing above. 

 For any cone $\tau \subset \mathbb{Z}^n \otimes \mathbb{R}$,
$\tau \cap M$ and $\tau \cap \mathbb{Z}^n$ are abelian semigroups.
We shall write $\mathbb{C}[\tau \cap M]$ and $\mathbb{C}[\tau \cap
\mathbb{Z}^n]$ for the $\mathbb{C}$-algebras generated by the corresponding 
Laurent monomials. Whenever we omit the lattice, writing
$\mathbb{C}[\tau]$, it should be assumed that the lattice is $M$. 

 The fan of $X$ in $L$ consists
of the single cone $L_+$, the dual of the cone $M_+$ of regular Laurent 
monomials in $M$ (similarly, we shall use $\mathbb{Z}^n_+$ and $\lbar_+$). 
The fan of any toric resolution of $X$ is given by 
a subdivision of $L_+$ into basic cones.

Fix such a toric resolution $Y$. Write $\mathfrak{F}$ for 
the set of basic cones which make up the fan
of $Y$. We shall denote by $A_\sigma$ the toric variety $\spec
\mathbb{C}[\sigma^\vee]$ corresponding to any cone $\sigma$ in
$L \otimes \mathbb{R}$. Then $Y$ is constructed in toric geometry
by gluing together $\{A_\sigma\}_{\sigma \in \mathfrak{F}}$: 
$A_{\sigma_1}$ and $A_{\sigma_2}$ are glued 
along $A_{\sigma_1 \cap \sigma_2} = \spec \mathbb{C}[(\sigma_1
\cap \sigma_2)^\vee]$.  Thus $\{ A_\sigma \}_{\sigma \in
\mathfrak{F}}$ is an open affine cover of $Y$.

  Now write $\mathfrak{E} \subset L$ for the set of all generators of 
these basic cones. In the toric geometry each element of
$\mathfrak{E}$ corresponds to either an exceptional divisor on $Y$ or
the proper transform of one of the coordinate hyperplanes in $X$. For
$e_i \in \mathfrak{E}$, write $E_i$ for the divisor on $Y$
corresponding to it.  

  It is often important whether the resolution is crepant or not.  
The discrepancy of each $E_i$ depends only on $e_i$ and not on the choice
of $Y$. If $e_i = (k_1, \dots, k_n) \in L$, then 
(\cite{ItReid96}, 1.4 and \cite{YPG87}, Prop. 4.8 for technicalities) 
the discrepancy of $E_i$ is $(\sum k_i)  -  1$, so the crepant divisors
correspond to the elements of $L$ which lie in the junior simplex:
$$ \Delta = \{ (k_1, \dots, k_n) \in L \otimes \mathbb{R}
\; |\; k_i > 0 \text{ and } \sum k_i = 1 \}$$  
Note that if a basic cone contains $e \in \Delta \cap L$, then $e$ must 
be one of its generators. So, for any resolution, $\Delta \cap L$ 
is a subset of $\mathfrak{E}$ and the crepant ones are precisely those for 
which this inclusion is an equality. 

\begin{exmpl}\label{exmpl-the-setup}

 Consider the group $G$ being $\frac{1}{8}(1,2,5)$, the group of 
$8$th roots of unity embedded into $\gsl_3(\mathbb{C})$ by 
$$\xi \mapsto 
\left( \begin{smallmatrix}
\xi^{1} &  &  \\ & \xi^{2} & \\ & & \xi^{5}
\end{smallmatrix} \right)
$$
We shall write $\chi_k$ for the character of $G$ given by $\xi \mapsto \xi^k$. 
So $x$ has weight $\chi_1$,  $y$ weight $\chi_2$ and $z$ weight $\chi_5$. 
 
 The lattice $L$ is generated in $(\mathbb{Z}^3)^\vee \otimes \mathbb{Q}$ 
by elements of $(\mathbb{Z}^3)^\vee$ and $\frac{1}{8}(1,2,5)$. 
The cone $L_+$, the positive octant, 
is the fan of $X$. A crepant resolution of $Y$ is given by 
a triangulation of the junior simplex $\Delta$ into basic
triangles. For the subsequent examples,  we choose the following
triangulation: 

\begin{center}
\includegraphics[scale=0.75]{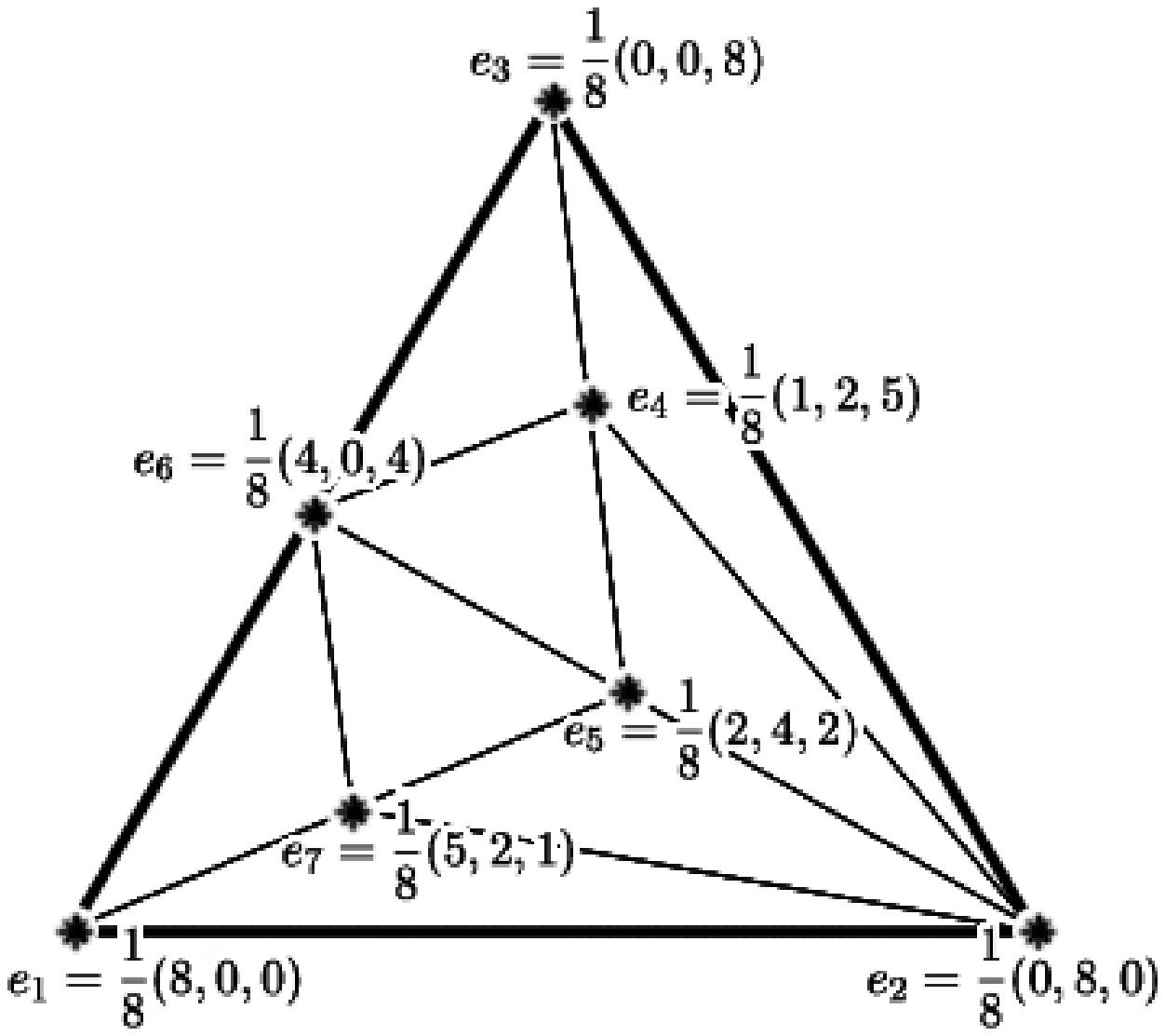}
\end{center}

So $\mathfrak{E} = \Delta \cap L = \{ e_1, \dots, e_7 \}$. 

And the basic cones of the fan $\mathfrak{F}$ of $Y$ are 
\begin{align*}
\mathfrak{F} = \biggl\{
& \left< e_1, e_2, e_7 \right>, \left< e_7, e_2, e_5 \right>, 
\left< e_4, e_2, e_5 \right>, \left< e_4, e_3, e_2 \right>, \\
& \left< e_3, e_4, e_6 \right>, \left< e_4, e_6, e_5 \right>, 
\left< e_6, e_5, e_7 \right>, \left< e_1, e_6, e_7 \right> \biggr\}
\end{align*}
This shall be the setup for all the subsequent examples.
\end{exmpl}

\subsection{Valuations} 

We now establish two simple results which translate the notions defined
in the Section \ref{section-valuations} into toric language. 

\begin{prps} \label{prps-toric-valuation}
 Let $Y$ be a toric resolution of $X$, $\mathfrak{F}$ its fan and
$\mathfrak{E}$ the set in $L$ of the generators of $\mathfrak{F}$. For any 
$e_i \in \mathfrak{E}$ and $m \in \mathbb{Z}^n$, 
\begin{align}
  v_{E_i}(x^m) = e_i(m) \in \mathbb{Q}
\end{align}
\end{prps}
\begin{proof}
 Take any basic cone $\sigma \in \mathfrak{F}$ such that 
$e_i \in \sigma$. Without loss of generality $i = 1$ and $\sigma =
\left< e_1, \dots, e_n \right>$. Let $\check{e}_1, \dots, \check{e}_{n}$ 
be the dual basis in $M$. 

 For any $m \in \mathbb{Z}^n$, $|G|m \in M$. Using the dual basis, 
$$ |G|m = \sum |G| e_i(m)\;\check{e}_i $$
 therefore 
$$ x^{|G|m} = (x^{\check{e}_1})^{|G|e_1(m)} \dots 
(x^{\check{e}_{n}})^{|G|e_{n}(m)} $$

 The restriction of the exceptional divisor $E_1$ to $A_\sigma$ is given
by the principal Weil divisor $(x^{\check{e_1}})$. Thus the local 
ring of $E_i$ is the coordinate ring of $A_\sigma$ localised at the
ideal $(x^{\check{e}_1})$, and so the valuation of $x^{|G|m} \in
\mathcal{O}_Y$ is $|G| e_1(m)$. By definition, 
$v_{E_1}(x^m) = \frac{1}{|G|} v_{E_1}(x^{|G|m}) = e_1(m)$. 
\end{proof}

The second result establishes which compatibility conditions a set of 
monomials $\{x^{m_\sigma}\}_{\sigma \in \mathfrak{F}}$ must satisfy for it 
to define a $G$-Cartier divisor. When the conditions are satisfied, 
we further establish the form which the corresponding $G$-Weil divisor must
take. 

\begin{prps} \label{prps-gens-of-weil}
A set $\{x^{m_\sigma}\}_{\sigma \in \mathfrak{F}} \subset
\mathbb{C}[\mathbb{Z}^n]$ of Laurent monomials defines a $G$-Cartier divisor 
$\{(A_\sigma, x^{m_\sigma})\}_{\sigma \in \mathfrak{F}}$ on $Y$ if and only if for any 
$e_i \in \mathfrak{E}$
\begin{align} \label{eqn-gluing}
e_i (m_\sigma) = e_i(m_{\tau}) \text{ for all } \sigma, \tau \ni e_i
\end{align}

When \eqref{eqn-gluing} holds, denote by $q_i$ the value
of $e_i(m_\sigma)$ for any $\sigma \ni e_i$. Then, under 
the isomorphism $\phi$ from Proposition \ref{prps-cartier-weil},
$\{(A_{\sigma},x^{m_\sigma})\}_{\sigma \in \mathfrak{F}}$ 
corresponds to the $G$-Weil divisor
$$ \sum_{e_i \in \mathfrak{E}} q_i E_i $$
\end{prps}
\begin{proof}
  Observe that if $\sigma, \tau \in \mathfrak{E}$ are such that 
$e_i$ belongs to both, then the generic point $p_{E_i}$ of $E_i$ 
lies in $A_{\sigma} \cap A_{\tau}$. If $\{(A_\sigma, x^{m_\sigma})\}$
is a $G$-Cartier divisor, then $x^{m_\sigma}/x^{m_{\tau}} \in 
\mathcal{O}^*(A_{\sigma} \cap A_{\tau})$, so we have
$v_{E_i}(x^{m_\sigma}/x^{m_{\tau}}) = 0$ and hence
$$ e_i(m_\sigma) = v_{E_i}(x^{m_\sigma}) = v_{E_i}(x^{m_{\tau}}) =
e_i(m_{\tau}) $$

  Conversely suppose we have $e_i(m_\sigma) = e_i(m_\tau)$ for all 
$e_i \in \sigma \cap \tau$. Then $m_\sigma - m_\tau \in (\sigma \cap 
\tau)^\perp$, and hence $x^{m_\sigma}/x^{m_\tau}$ is invertible 
in $\mathbb{C}[(\sigma \cap \tau)^\vee] = 
\mathcal{O}_Y(A_\sigma \cap A_\tau)$
as required.

  For the last part, recall that $\phi(\{(A_\sigma, x^{m_\sigma})\})$
is defined as the sum $\sum n_D D$ over all prime divisors on $Y$ where $n_D =
v_D(x^{m_\sigma})$ for any $\sigma$ such that $D \cap A_\sigma \neq
\emptyset$. So it suffices to prove that, for all $\sigma \in
\mathfrak{F}$, the restrictions of the
principal divisor $(x^{m_\sigma})$ and $\sum_{i \in \mathfrak{E}} q_i
E_i$ to $A_\sigma$ are identical. 

  Without loss of generality, we can take
$\sigma = \left< e_1, \dots, e_n \right>$. Then $\mathcal{O}_{A_\sigma} =
\mathbb{C}[t_1, \dots, t_n]$ where $t_i = x^{\check{e}_i}$. 
We have $x^{m_\sigma} = \prod_{e_i \in \sigma} t_i^{q_i}$ and 
recall (proof of Proposition \ref{prps-toric-valuation}) 
that $E_i|_{A_\sigma} = (t_i)$. Therefore 
$$(x^{m_\sigma})|_{A_\sigma} = \sum_{e_i \in \sigma} q_i\; (t_i) =  
(\sum_{e_i \in \sigma} q_i \; E_i)|_{A_\sigma} $$      
and the result follows.
\end{proof}
\begin{remarks*}
\begin{enumerate}
\item Observe that the `only if' part of the proof is completely
general and doesn't rely on the toric technology. It is the standard 
argument used to show that the morphism $\phi$ taking Cartier
divisors to Weil divisors is well-defined.  

On the other hand the `if' argument is toric-specific and relies heavily
on the fact that the invertible functions on $A_\sigma \cap A_\tau$
are precisely the monomials in $(\sigma \cap \tau)^\vee$.

\item Note that, in particular, we have proved that for any $m \in
\mathbb{Z}^n$, the sum 
$$ \sum_{i \in \mathfrak{E}} v(E_i, x^{m}) E_i $$
 is a valid $G$-Weil divisor on $Y$. 
Recalling the definition of $G$-Weil divisors, 
this provides an independent proof that for any prime divisor $D$ which is 
not $E_i$ for some $i \in \mathfrak{E}$, we have $$v(D,\chi) = 0$$ 
for all $\chi \in G^\vee$, since $v(D,\chi)$ is defined as the
fractional part of the valuation of any homogeneous rational 
function of weight $\chi$ on $D$. 
\end{enumerate}
\end{remarks*}

\begin{exmpl} \label{exmpl-weil-divisor}
To illustrate the above, in the context of the Example \ref{exmpl-the-setup}, we shall calculate explicitly the $\chi_6$-Cartier
divisor corresponding to the $\chi_6$-Weil divisor
$$ D =  \frac{7}{4}E_4 + \frac{1}{2}E_5 - \frac{1}{4} E_7 $$

 Consider the cone $\sigma = \left< e_4, e_5, e_6 \right>$. Calculating 
the dual basis which generates the abelian semigroup $\check{\sigma} \cap M$,
we get 
$$ \check{e}_4 = (-2, 0, 2), 
\quad \check{e}_5 = (1, 2, -1), \quad \check{e}_6 = (2, -1, 0) $$

So $A_\sigma = \spec \mathbb{C}[\frac{z^2}{x^2}, \frac{x y^2}{z},
\frac{x^2}{y}]$ and the restrictions of $E_4$, $E_5$ and $E_6$ to 
$A_\sigma$ are given by $(\frac{z^2}{x^2})$, $(\frac{x y^2}{x^2})$ 
and $(\frac{x^2}{y})$ respectively. To specify $D$ on $A_\sigma$ we
need $f \in K{\chi_6}(\mathbb{C}^3)$ such that 
$v_{E_4}(f) = \frac{7}{4}$, $v_{E_5}(f) = \frac{1}{2}$ and $v_{E_6}(f)
= 0$, so we take
$$\left(\frac{z^2}{x^2}\right)^{7/4} \left(\frac{x
y^2}{z}\right)^{1/2} \left(\frac{x^2}{y}\right)^0 =
\frac{z^3 y }{x^3} $$
to be $f$. 

  Repeating the same calculations for the remaining cones in the fan 
$\mathfrak{F}$ we get the $\chi_6$-Cartier divisor given by 
\begin{center}
\includegraphics[scale=0.50]{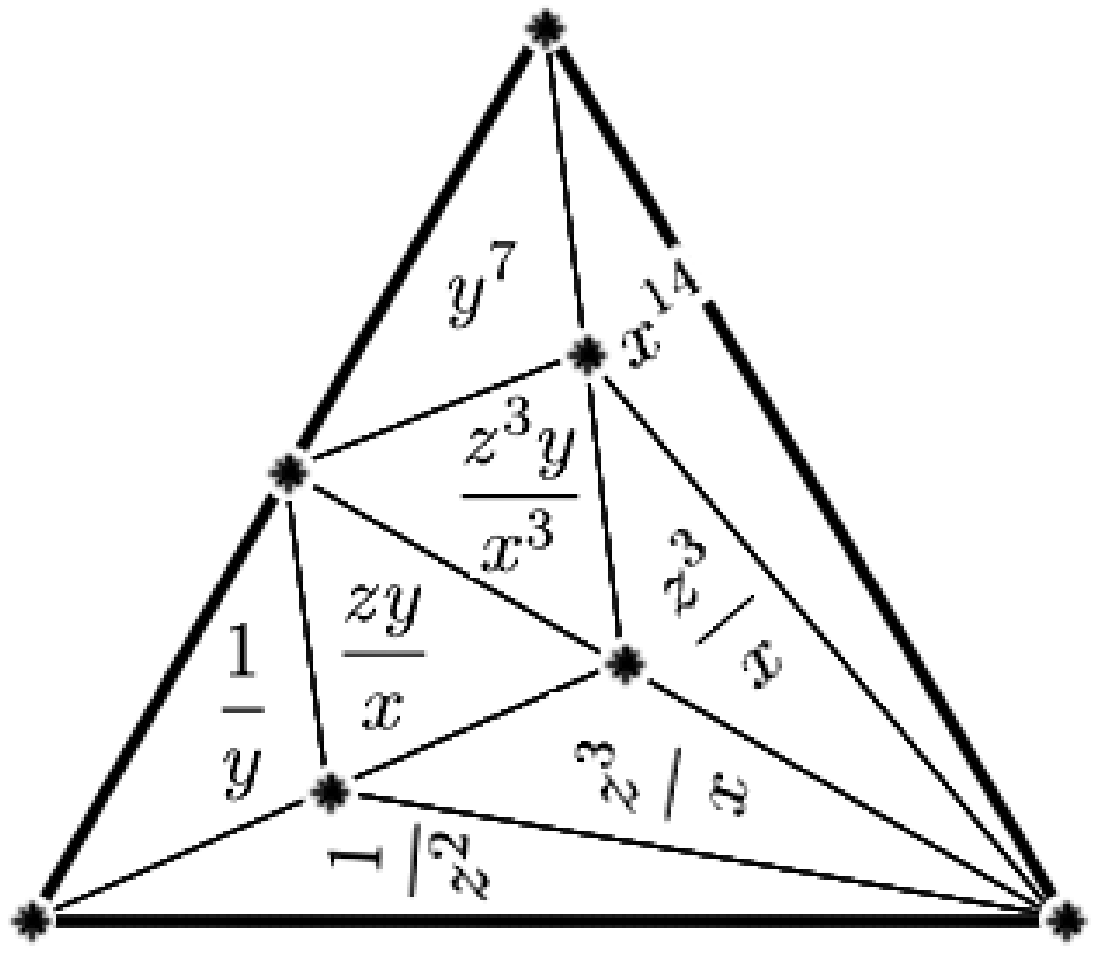}
\end{center}
and we can indeed see that, as all the monomials representing the
divisor have weight $\chi_6$, their ratios are all invariant
and the sub-$\mathcal{O}_Y$ module of $K(\mathbb{C}^n)$ they generate
is an invertible sheaf on $Y$ with the natural action of $G$ by $\chi_2$.
\end{exmpl}

\subsection{Representations of the McKay Quiver} \label{section-quivers}
  We now introduce a useful way to visualise the mechanics of 
a family of $G$-constellations over a particular toric affine piece of $Y$. 
Suppose we have a family $\mathcal{F}$ of $G$-constellations on $Y$
and a cone $\sigma$ in the fan $\mathfrak{F}$. In this section, 
we are interested in looking up close at the structure of $\mathcal{F}$ 
restricted to the corresponding affine piece $A_\sigma$. 

  Over $A_\sigma$ the sheaf $\mathcal{F}$ is trivialised and we have
$$ \mathcal{F}(A_\sigma) \simeq \mathbb{C}[\sigma^\vee]
\otimes_\mathbb{C}
\regrep \simeq \bigoplus_{\chi} F_\chi $$
where each $F_\chi$ is isomorphic to $\mathbb{C}[\sigma^\vee]$ and $G$ acts 
on it by $\chi$.
Evidently, the whole structure of $\mathcal{F}$ as a family of
$G$-constellations on $A_\sigma$ is contained in the way that $\regring$ acts 
on $F_\chi$s. An effective method to visualise the mechanics of this is 
to consider the representations of \it the McKay quiver of $G$\rm. 
We shall briefly summarize the necessary background. For a more 
detailed exposition of the following material see [  ]. 

\begin{defn}
A \tt quiver \rm consists of a vertex set $Q_0$, 
an arrow set $Q_1$ and two maps $h: Q_1 \rightarrow Q_0$ and 
$t: Q_1 \rightarrow Q_0$ giving the head $hq \in Q_0$ and the tail $tq
\in Q_0$ of each arrow $q \in Q_1$. 
\end{defn}

\begin{defn}
Let $G$ be a finite subgroup of $\gl (\givrep)$. 
Then the \tt McKay quiver of $G$ \rm is 
the quiver with the vertex set $Q_0$ labelled by the irreducible
representations $\rho$ of $G$ and the arrow set $Q_1$ which 
has precisely $\dim \homm_G (\rho_i, \rho_j \otimes \givrep)$ arrows
going from the vertex $\rho_i$ to the vertex $\rho_j$. 
\end{defn}

\begin{exmpl}
\begin{enumerate}
\item In our case of $G$ being abelian and $\givrep$ identified with 
$\mathbb{C}^n$, we have a decomposition of $\givrep^\vee$ into
irreducible representations as $\oplus \mathbb{C} x_i$,
where $x_i$s are the basic monomials. Then, writing $U_\chi$ for 
the  representation corresponding to $\chi \in G^\vee$
\begin{align}
\homm_G(U_{\chi_i}, U_{\chi_j} \otimes \mathbb{C}^n) = 
\bigoplus_{x_k \;|\; \chi_i \rho^{-1}(x_k) = \chi_j} 
\homm_G(x_k \otimes U_{\chi_i}, U_{\chi_j})
\end{align}
where by $x_k \otimes U_{\chi_i}$, we denote the space 
$\mathbb{C}x_k \otimes_\mathbb{C} U_{\chi_i}$. Each of the spaces 
$\homm_G(x_k \otimes U_{\chi_i}, U_{\chi_j})$ is one-dimensional and
so has one arrow from $\chi_i$ to $\chi_j$ corresponding to it. Thus 
the quiver consists of $|G|$ vertices labelled by characters $\chi \in
G^\vee$ and out of each vertex $\chi$ emerge $n$ arrows, each
corresponding to one of the one-dimensional spaces 
$\homm_G(x_k \otimes U_\chi, U_{\chi \rho(x_k)})$. We shall 
write $(\chi, x_k) \in Q_1$ to denote such an arrow. 

\item For a concrete example, the reader can verify that the McKay
quiver for $G = \frac{1}{8}(1,2,5)$ (see Example
\ref{exmpl-the-setup}) looks like:

\begin{center}
\includegraphics[scale=0.45]{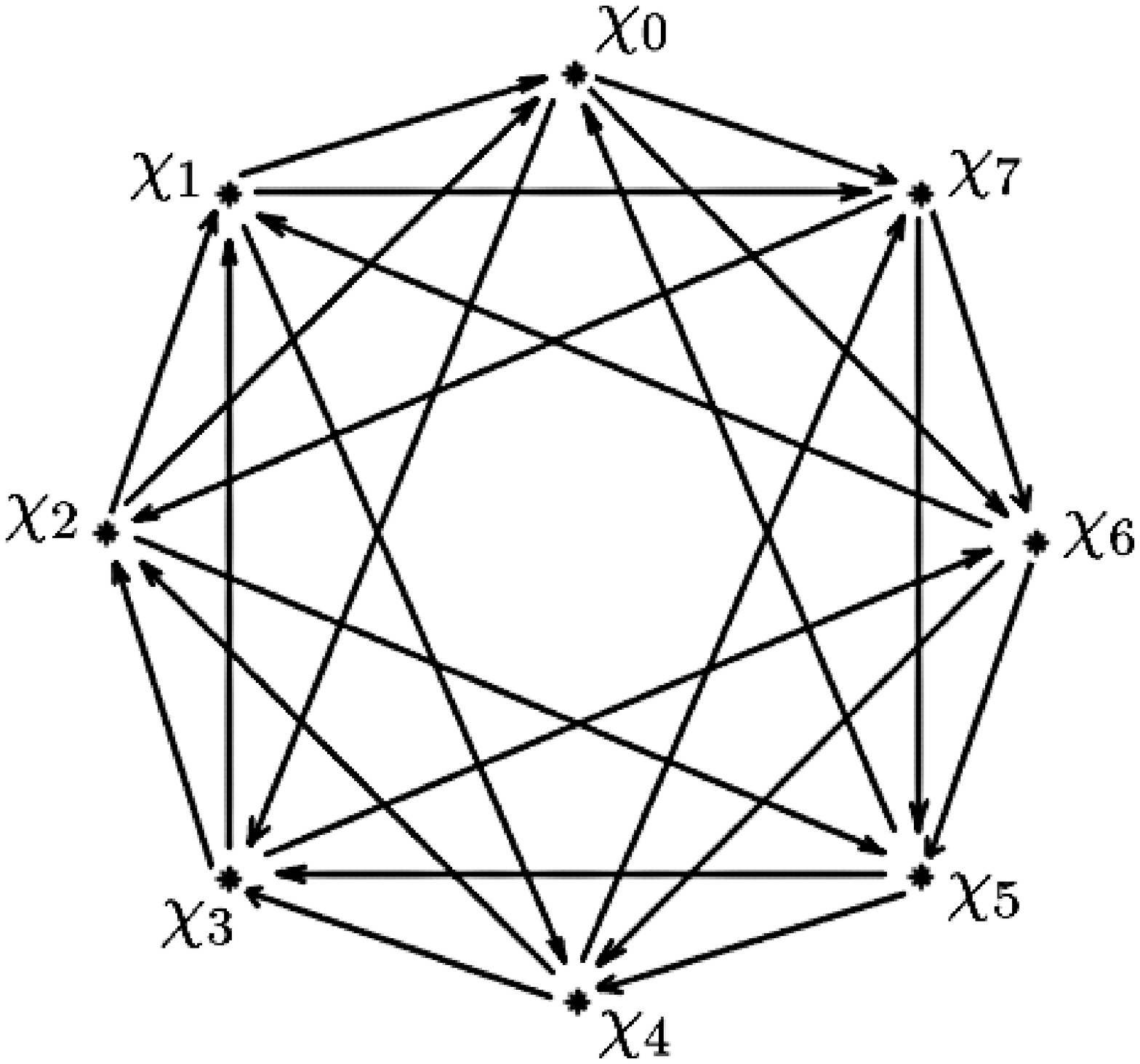}
\end{center}
  
\end{enumerate}
\end{exmpl}

  A good reason for contemplating the McKay quiver of $G$ is that it is 
possible to establish a 1-to-1 correspondence between a subset of its \it
representations \rm and $(G,\regring)$-modules.

\begin{defn}
A \tt representation \rm of a quiver is a graded vector space
$\oplus_{i \in Q_0} V_i$ and a collection 
$\{\alpha_q:\; V_{tq} \rightarrow V_{hq}\}_{q \in Q_1}$ of
linear maps indexed by the arrow set of the quiver. 
A morphism from  $(\oplus V_i, \{ \alpha_q \} )$ to 
$( \oplus V'_i, \{ \alpha'_q \} )$ 
is a collection of linear maps 
$\{\theta_i : V_i \rightarrow V'_i\}_{i \in
Q_0}$ forming commutative squares with $\alpha_q$s and $\alpha'_q$s.
\end{defn}

  Given a $G$-representation $V$, it is traditional, in case of 
$G$ being a general finite subgroup of $\gl_n$, to consider representations 
of the McKay quiver on a graded vector space $\oplus V_\rho$ where
$V_\rho = \homm_G(\rho, V)$. It is then possible (\cite{Sar-In96b}) 
to establish a 1-to-1 correspondence between such representations and 
elements of $\homm_G (\givrep^\vee \otimes V, V)$. And, in the light of 
the remarks after the Definition $\ref{defn-gcon}$, there is a 1-to-1 
correspondence between all the $(G,\regring)$-module structures on $V$ and 
the elements of $\homm_G(\givrep^\vee \otimes V, V)$ which satisfy the
commutator relations \eqref{eqn-commutator-conditions-gen}.

  However, in the case when the group $G$ is abelian, a considerable 
shortcut can be taken by considering the representations directly 
into graded vector space $\oplus V_\chi$, where $V_\chi$ is 
the $\chi$-eigenspace of $V$. 
We again have the correspondence between representations of McKay quiver 
on $\oplus V_\chi$ and elements of $\homm_G(\givrep^\vee \otimes V,
V)$ and consequently the correspondence with $G$-constellations. 
Explicitly, if we have a $(G,R)$-structure on $V$, then the action 
map $V \rightarrow V$ for each basic monomial $x_i$ is $G$-equivariant
and so splits into maps $V_\chi \rightarrow V_{\chi / \rho(x_i)}$. Each 
such map gives precisely the map $\alpha_{\chi,x_i} \in 
\homm(V_\chi, V_{\chi / \rho(x_i)})$ in the corresponding 
representation of the quiver.

 In case of $V = \regrep$, if we make an explicit choice of 
a basis vector $e_\chi$ for each $V_\chi$, this gives us bases for all
$\homm_G(x_i \otimes V_\chi, V_{\chi / \rho(x_i)})$. Then every
McKay quiver representation on $\oplus V_\chi$ gains a unique map $\xi: Q_1
\rightarrow \mathbb{C}$ associated with it, defined by 
\begin{align*}
\alpha_{\chi, x_i}(e_\chi) = \xi(\chi, x_i) e_{\chi / \rho(x_i)}
\end{align*}

Considering a family of $G$-constellations $\mathcal{F}$ parametrised by 
an affine piece $A_\sigma$ of $Y$, we have, as outlined in the beginning of the 
section, 
$$\mathcal{F}(A_\sigma) \simeq 
\mathbb{C}[\sigma^{\vee}] \otimes_\mathbb{C} \regrep $$ 
We then write the $\chi$-eigenspace decomposition 
$\mathcal{F}(A_\sigma) = \oplus F_\chi$, and all the correspondences above 
work just as well with $\mathbb{C}[\sigma^\vee]$-modules as they did 
with complex vector spaces.

This technology presents us with a compact way to write down
the $\regring$-module structure on $\mathcal{F}|_{A_\sigma}$. 
After a choice of bases, a representation of the McKay quiver 
becomes a map $\xi: Q_1 \rightarrow \mathbb{C}[\sigma^\vee]$ readily 
pictured as a McKay quiver of $G$ with $\xi(\chi, x_i)$ written 
above each arrow $(\chi, x_i) \in Q_1$.  
In this way it is also easy to calculate explicitly the
$G$-constellation in $\mathcal{F}$ parametrised by any point of
$A_\sigma$. If a point $p \in A_\sigma$ is defined by a map 
$\ev_p: \mathbb{C}[\sigma^\vee] \rightarrow \mathbb{C}$, 
then the corresponding quiver representation is given by the map 
$\xi_p = \ev_p \circ\; \xi: Q_1 \rightarrow \mathbb{C}$.

Finally, let us consider deformations of the generic orbit. If $\mathcal{F}$
is one such, then it comes with an embedding $\iota: \mathcal{F}
\rightarrow K(\mathbb{C}^n)$. 
Its image $\iota(\mathcal{F})$ splits into $\chi$-eigenspaces, 
which are invertible sheaves, so we can take a set 
$\{f_\chi\} \in K(\mathbb{C}^n)$,
where each $f_\chi$ is homogeneous of weight $\chi$ and a genereator 
of $\chi^{-1}$-eigenspace of $\mathcal{F}$ over $A_\sigma$. 
The $\regring$-module structure comes for free with the embedding into
$K(\mathbb{C}^n)$ and the corresponding quiver representation is
given by the map 
$\xi: Q_1 \rightarrow \mathbb{C}[\sigma^\vee]$ defined by 
$$ (\chi^{-1}, x_i ) \mapsto \frac{x_i f_\chi}{f_{\rho(x_i)\chi}} $$  
with respect to the choice of generators $f_\chi$.
  
\begin{exmpl}
 Let us work through an actual example. Let $G =
\frac{1}{8}(1,2,5)$ and $\sigma = \left< e_4, e_5, e_6 \right>$.
Recall from the Example \ref{exmpl-weil-divisor} that the calculation of
the dual basis in $M$ gives us the local coordinates on $A_\sigma = \spec
\mathbb{C}[\sigma^\vee]$ as $\mathbb{C}[\sigma^\vee] = \mathbb{C}[\frac{z^2}{x^2}, \frac{x y^2}{z},
\frac{x^2}{y}]$.
  
Consider $\mathcal{F} = \oplus_{\chi_i \in G^\vee} \mathcal{O}_{A_\sigma} f_i \subset
K(\mathbb{C}^n)$ where
\begin{align*}
f_0 = 1 && f_1 = x && f_2 = y \\
f_3 = xy  && f_4 = \frac{z}{x} && f_5 = z \\
f_6 = \frac{yz}{x} && f_7 = yz &&
\end{align*}

 Now for any choice of $f_i$, as long as each $f_i \in
K_{\chi_i}(\mathbb{C}^n)$, the generic fiber $\oplus K(Y) f_i$
is the whole of $K(\mathbb{C}^n)$. The latter has a natural structure
of a $G$-constellation, and so it has a corresponding 
quiver representation. Let $\xi': Q_1 \rightarrow K(Y)$ be the map specifying 
it with respect to $\{f_i\}$s as the choice of eigenspace bases. 

 We claim that $\mathcal{F}$ is closed under $R$-action in $K(\mathbb{C}^n)$ 
and hence defines a family of $G$-constellations parametrised by $A_\sigma$. 
We shall verify this statement in the course of calculating the map
$\xi'$ and seeing that it restricts to a map $Q_1 \rightarrow
\mathbb{C}[\sigma^\vee]$, which defines the quiver representation
corresponding to our family. 
 
 Consider the arrow $(\trch_, x)$. As described above, in the corresponding 
quiver representation the map
$K(Y) f_{0} \rightarrow  K(Y) f_{1}$
is given by multiplication by $x$. Hence we get $$f_0 \mapsto
1\; f_1$$ and so we label this arrow by $$1 =
\left(\frac{z^2}{x^2}\right)^0 \left(\frac{x y^2}{z}\right)^0 
\left(\frac{x^2}{y}\right)^0$$ Similarly 
the arrow $(\chi_5, z)$ corresponds to the map $f_3 \mapsto
xyz\; f_0$ and so we label it by $$ xyz =
\left(\frac{z^2}{x^2}\right)^1 \left(\frac{x y^2}{z}\right)^1 
\left(\frac{x^2}{y}\right)^1 $$

 Repeating this for all the arrows of the quiver we obtain:
\begin{center}
\includegraphics[scale=0.32]{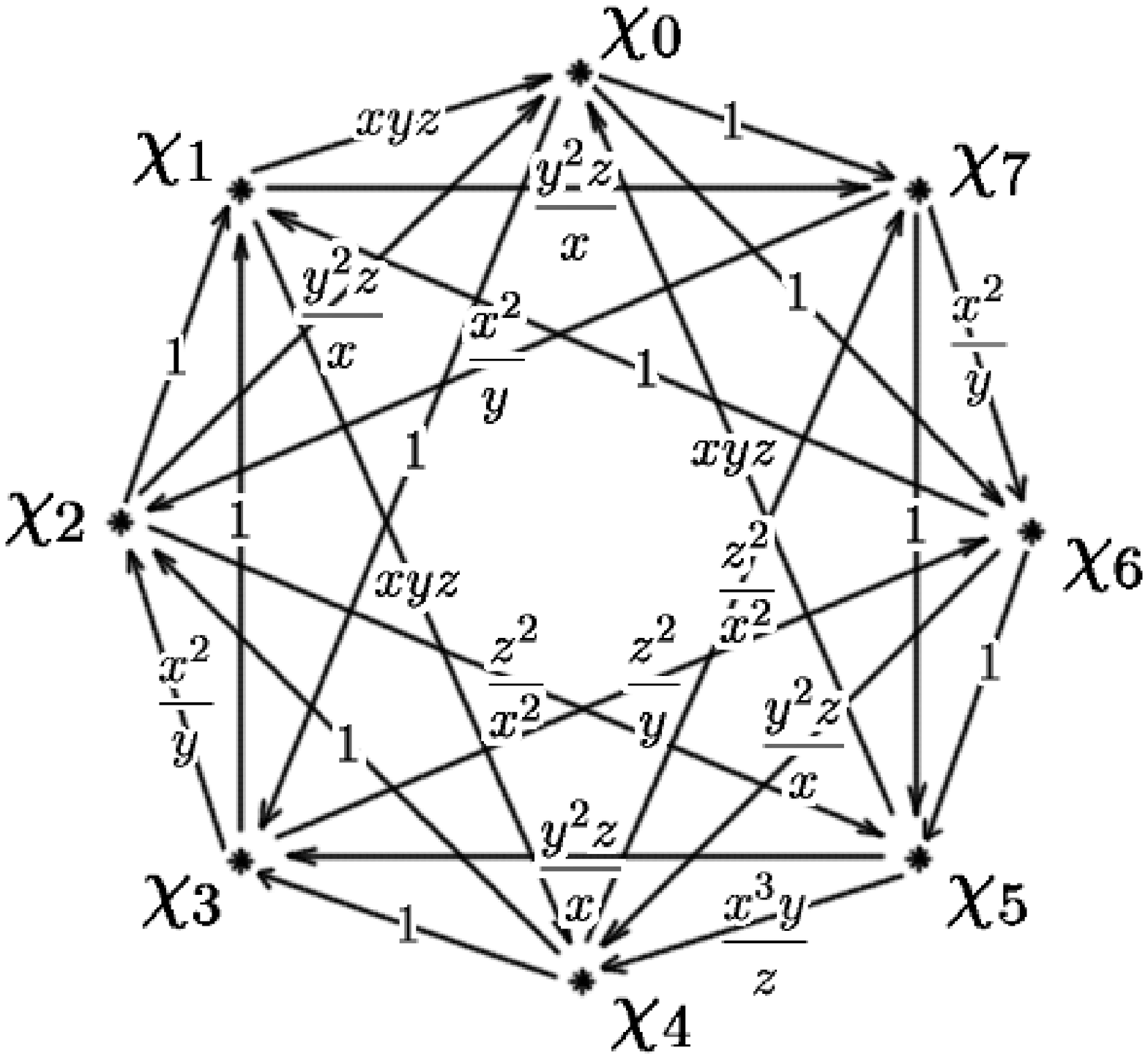}
\includegraphics[scale=0.32]{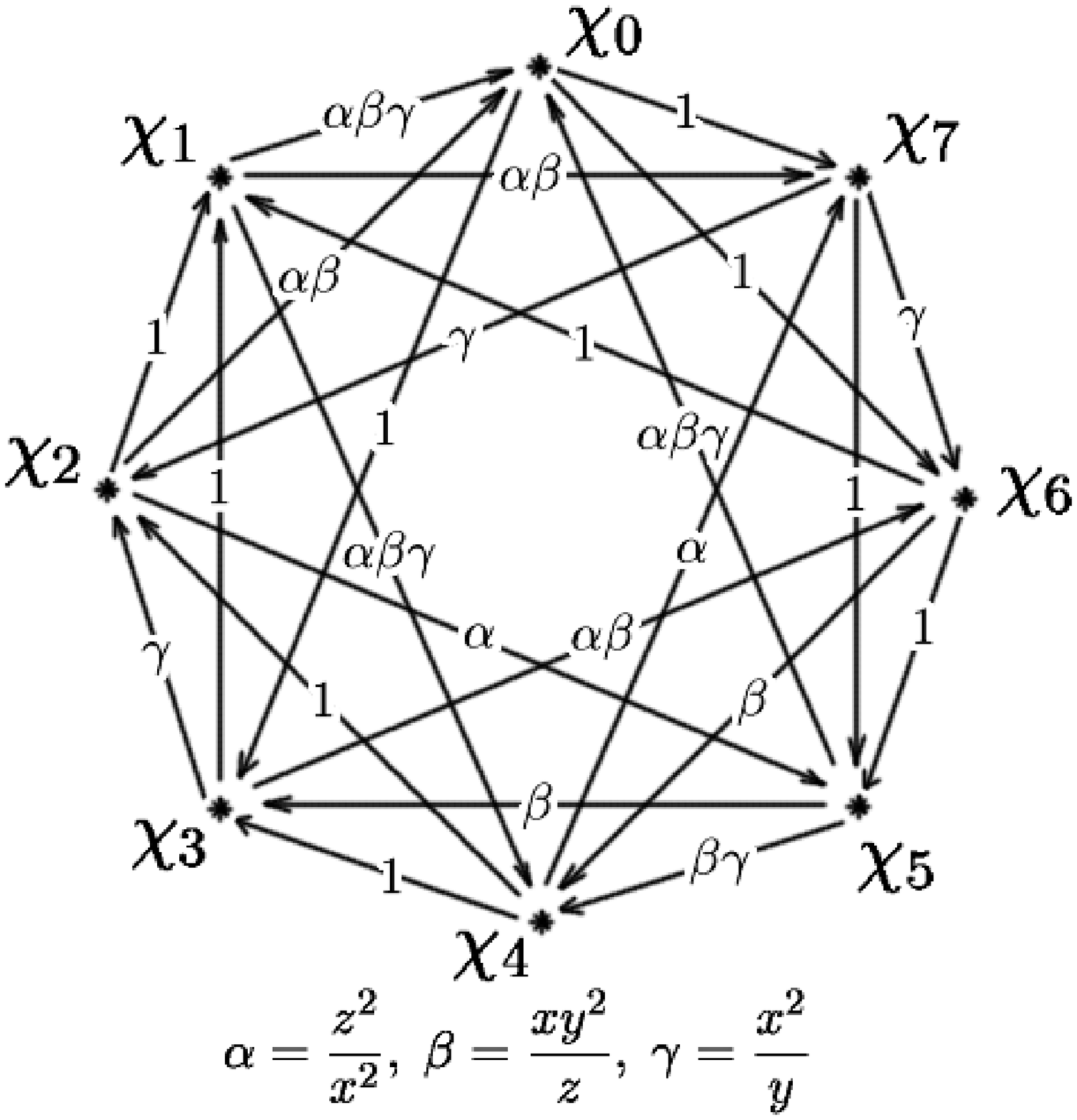}
\end{center}

 In the diagram on the right we have written all the functions
marking the arrows in terms of positive powers of the local coordinates 
$\alpha, \beta, \gamma$
on $A_\sigma$. This demonstrates that we indeed have a map 
$$ \xi: \quad Q_1 \rightarrow 
 \mathbb{C}[\sigma^\vee] \quad = \mathbb{C}\left[\frac{z^2}{x^2},
\frac{x y^2}{z}, \frac{x^2}{y}\right]$$ 
so $\mathcal{F}$ is indeed a family of $G$-constellations parametrised
by $A_\sigma = \spec[\alpha, \beta, \gamma]$. The $G$-constellations
parametrised by each point of $A_\sigma$ are readily calculated by assigning 
specific values to $\alpha$, $\beta$ and $\gamma$ in the diagram 
on the right.

\end{exmpl}

\section{Reductors} \label{section-reductors}

\subsection{Reductor Pieces} \label{subsection-reductors}

 As in Section \ref{section-quivers}, let $Y$ be a toric 
resolution, $\sigma \in \mathfrak{F}$ a cone in its fan and
$\mathcal{F}$ a deformation of the generic orbit across $Y$. 
If we have a set of generators $\{f_\chi \;|\; f_\chi \in \kgchif_ \}$ 
such that 
$$ \iota(\mathcal{F})(A_\sigma) = \bigoplus \mathbb{C}[\sigma^\vee] f_\chi $$
then we must have 
\begin{align} \label{eqn-red-piece-cond}
\frac{x_i f_\chi}{f_{\rho(x_i) \chi}} \in \mathbb{C}[\sigma^\vee] 
\end{align}
for all basic monomials $x_i$ and $\chi \in G^\vee$. 

  But observe that, conversely, for any set $\{f_\chi \;|\; f_\chi \in
\kgchif_\}$ for which $\eqref{eqn-red-piece-cond}$ holds, 
the $\mathbb{C}[\sigma^\vee]$-submodule of $K(\mathbb{C}^n)$
generated by $f_\chi$ is closed under the natural action of 
$\regring$ on $K(\mathbb{C}^n)$ by multiplication. It is certainly
closed under the $G$-action, so it is a $(G,\regring)$-submodule of
$K(\mathbb{C}^n)$ and a family of $G$-constellations 
parametrised by $A_\sigma$. 

  This observation motivates the rest of this section. But first
we make a useful definition
\begin{defn}
A \tt reductor piece for a basic cone $\sigma
\subset L$ \rm of the fan $\mathfrak{F}$ of the toric resolution $Y$ is 
a set $\{ f_\chi \;|\;  f_\chi \in \kgchif_ \}$ 
such that for any basic monomial $x_i$ and any $\chi \in G^\vee$
we have 
\begin{align}
\frac{x_i f_{\chi} }{f_{\rho(x_i)\chi}} \in \mathbb{C}[\sigma^{\vee}]
\end{align}
\end{defn}

Thus, if we wanted to explicitly construct a family of
$G$-constellations parametrised by $Y$, we could do it by 
producing a reductor piece for each cone $\sigma$ in the fan
$\mathfrak{F}$. Every such would give a family of $G$-constellations 
parametrised by open affine piece $A_\sigma$. However, we would 
need these families to `glue together', i.e. the restrictions to 
$A_\sigma \cap A_{\sigma'}$ of the families generated on 
$A_\sigma$ and $A_{\sigma'}$, respectively, must be isomorphic 
for any two cones $\sigma, \sigma' \in \mathfrak{F}$. 
The general way to guarantee this is independent of the toric technology 
altogether, taking us back to $G$-Weil divisors and to where 
Section \ref{section-valuations} left off.  

\subsection{Reductor Sets}

 From now on $Y$ is once again an arbitrary, not necessarily toric, 
resolution of $X$.

  Let $\mathcal{F}$ be a deformation of the generic orbit. It comes 
with a choice of an embedding $\iota':~\mathcal{F} \hookrightarrow 
K(\mathbb{C}^n)$.  Then $\mathcal{F}$ splits into $G$-eigensheaves 
as $\oplus \mathcal{F}_\chi$ and, as
per Section \ref{section-valuations}, each $\mathcal{F}_\chi$ 
defines a linear equivalence class of $\chi$-divisors embedding it into 
$K(\mathbb{C}^n)$, and $\iota'(\mathcal{F}_\chi)$ pinpoints a specific
element of that class. Hence  
$\iota'(\mathcal{F}) = \oplus_\chi \mathcal{L}(-D_\chi)$
for some unique set of $G$-divisors $\{D_\chi\}_{\chi \in G^\vee}$.
Note that it is important here that $\mathcal{L}(-D_\chi)$ is not
merely an abstract line bundle corresponding to $-D_\chi$, but 
a specific sub-$\mathcal{O}_Y$-module of $K(\mathbb{C}^n)$ as 
per its definition. 

Thus each subsheaf of the constant sheaf $K(\mathbb{C}^n)$ on $Y$, which 
is an image of an isomorphism class of deformations of the generic orbit, 
is of the form $\oplus \mathcal{L}(-D_\chi)$, where each $D_\chi$ 
is a $\chi$-divisor on $Y$. 

\begin{lemma} \label{lemma-iso-dofgs}
Let $\mathcal{F} = \oplus \mathcal{L}(-D_\chi)$ and $\mathcal{F}' = 
\oplus \mathcal{L}(-D'_\chi)$ be two deformations of the generic orbit 
across $Y$. Then they are isomorphic as sheaves of
$(G,\regring)$-modules if and only if there exists $g \in K(Y)$ such that 
\begin{align}
D'_\chi - D_\chi = (g)
\end{align}
for all $\chi \in G^\vee$. 
\end{lemma}
\begin{proof}
The `if' part is immediate, observe that we have a natural isomorphism
$\mathcal{L}(A) \otimes \mathcal{L}(B) \rightarrow
\mathcal{L}(A+B)$ given by multiplication in $K(\mathbb{C}^n)$. 
Applying this to $ - D_\chi - (g) = - D'_\chi$ yields isomorphism 
$\mathcal{F} \rightarrow \mathcal{F'}$ given by $s \mapsto s/g$. 

For the `only if' part, 
let $\phi: \oplus \mathcal{L}( - D_\chi) \rightarrow \oplus
\mathcal{L}(- D'_\chi)$ be a $(G,\regring)$-equivariant isomorphism.
Then it restricts to $\phi_{\chi}: \mathcal{L}( - D_\chi)
\xrightarrow{\sim} \mathcal{L}( - D'_\chi)$ for all $\chi \in G^\vee$.
Then $\phi_{\chi}$ induces a map $\mathcal{L}(0) \xrightarrow{\sim}
\mathcal{L}( - D'_\chi + D_\chi)$, so let $g_{\chi} \in K(\mathbb{C}^n)^G$ 
be an image of $1$ under this map. 
Then $D'_\chi - D_\chi = (g_{\chi})$ and $\phi_{\chi}$ is 
given by $s \mapsto g_{\chi} s$ for any $s \in \mathcal{L}( - D_\chi)$.

It remains to show that all $g_\chi$ are equal. Fix any $\chi \in G^\vee$ 
and consider any $G$-homogeneous $m \in R$ of weight $\chi$.
Take any $s \in \mathcal{L}( - D_\trch_) \subset K(\mathbb{C}^n)$.  
Then $m s \in \mathcal{L}(- D_\chi)$ and using $\regring$-equivariance
of $\phi$
\begin{align}
\phi(ms) = m \phi(s) = g_{\trch_} m s
\end{align}
and hence $g_{\chi} = g_{\trch_}$ for all $\chi \in G^\vee$. 
\end{proof}

\begin{cor} \label{cor-equiv-families-divisor-diff}
Let $\mathcal{F} = \oplus \mathcal{L}( - D_\chi)$ and $\mathcal{F}' = 
\oplus \mathcal{L}( - D'_\chi)$ be two deformations of the generic orbit 
across $Y$. Then they are equivalent if and only if there 
exists a $\trch_$-divisor $N$ such that 
\begin{align}
D'_\chi - D_\chi = N
\end{align}
for all $\chi \in G^\vee$. 

\end{cor}

\begin{proof}
Once again, the `if' direction is immediate: an isomorphism 
$\mathcal{F} \otimes \mathcal{L}( - N) \rightarrow \mathcal{F}'$ is
given by multiplication in $K(\mathbb{C}^n)$.  

Conversely, if the families are equivalent then let 
$\mathcal{N}$ be an invertible sheaf on $Y$ such that 
$\mathcal{F}' \simeq \mathcal{F} \otimes \mathcal{N}$. Choose any Weil 
divisor $N'$ such that $\mathcal{N} = \mathcal{L}(- N')$. Then apply Lemma
\ref{lemma-iso-dofgs} to the isomorphic families 
$\oplus \mathcal{L}( - D_\chi - N')$ and $\mathcal{L}( - D'_\chi)$ to 
obtain $g \in K(\mathbb{C}^n)$ such that $D'_\chi - D_\chi - N' = (g)$
for all $\chi \in G^\vee$. Setting $N = N' + (g)$ finishes the proof.
\end{proof}

\begin{cor} \label{cor-equiv-classes-normalized-sets}
In every equivalence class of deformations of the generic orbit there
exists a unique family $\mathcal{F}$ of the form 
$\oplus \mathcal{L}( - D_\chi)$ with $D_{\trch_} = 0$.
\end{cor}

\begin{proof}
 Given an arbitrary deformation of the generic orbit $\mathcal{F}$ we 
can find an isomorphic family of the form $\oplus \mathcal{L}(-D_\chi)$.
Then setting $D'_{\chi} = D_{\chi} - D_{\trch_}$ we obtain an
equivalent family $\mathcal{L}(-D'_{\chi})$ with the required properties. 
Finally, Corollary \ref{cor-equiv-families-divisor-diff} shows the uniqueness. 
\end{proof}

In the view of all of the above, we make following definitions:

\begin{defn}
Let $\{D_\chi\}_{\chi \in G^\vee}$ be a set of $G$-divisors. We shall 
call it a \tt prereductor set \rm if each $D_\chi$ is a $\chi$-Weil 
divisor. We shall call it a \tt reductor set \rm if 
$\oplus \mathcal{L}( - D_\chi)$ with the inclusion map into
$K(\mathbb{C}^n)$ is a deformation of the generic orbit. We shall 
say the reductor set is \tt normalised \rm if $D_{\trch_} = 0$.  
\end{defn}

\subsection{Reductor Condition}

 We have seen that a deformation of the generic orbit can be specified 
(up to an isomorphism) by a set of $G$-Weil divisors on $Y$ which gives its 
embedding into $K(\mathbb{C}^n)$. Here we investigate an opposite question: 
for which prereductor sets $\{D_\chi\}$ is $\oplus \mathcal{L}(-D_\chi)$  
a family of $G$-constellations.

  We observe that $\oplus \mathcal{L}( - D_\chi)$ is always a 
sub-$\mathcal{O}_Y$-module of $K(\mathbb{C}^n)$ closed under the
$G$-action. However, for a general choice of divisors
$D_\chi$, there is no guarantee that the $\oplus \mathcal{L}( - D_\chi)$
will be closed under the $\regring$-action on $K(\mathbb{C}^n)$. 

\begin{prps}[Reductor Condition] \label{prps-reductor-condition}
Let $\{D_\chi\}$ be a prereductor set. Then it is a reductor set if 
and only if, for any $f \in \regring_G$, a G-homogeneous polynomial, 
the divisor 
\begin{align} \label{eqn-reductor-condition} 
D_\chi + (f) - D_{\chi \rho(f)} \geq 0
\end{align}
i.e. it is effective. 
\end{prps}

\bf Remarks: \rm 
\begin{enumerate}
\item It is, of course, sufficient to check 
\eqref{eqn-reductor-condition} only for $f$ being one of the 
basic monomials $x_1, \dots, x_n$. This leaves us 
with a finite number of inequalities to check. Note also that 
the principal divisor $(x_j)$ is very easy to compute in toric
case. It follows immediately from Proposition \ref{prps-gens-of-weil}
that it is $\sum_{e_i \in \mathfrak{E}} e_i(x_j) E_i$. Observe that 
$e_i(x_j)$ is simply the $j$th coordinate of $e_i$ in $L$. 

\item Numerically, if we write each $D_\chi$ as $\sum q_{\chi,P} P$,
each inequality $\eqref{eqn-reductor-condition}$ becomes a set of
inequalities \begin{align} \label{eqn-red-cond-num} q_{\chi,P} +
v_P(f) - q_{\chi \rho(f), P} \geq 0 \end{align} for all prime divisors
$P$ on $Y$.  The important thing to notice here is that the subsets of
inequalities for each prime divisor $P$ are all independent of each
other. We can speak of $\{D_\chi\}$ satisfying or not satisfying the
reductor condition at a given prime divisor $P$. Moreover, we can 
construct reductor sets $\{D_\chi\}$ by independently choosing for each prime
divisor $P$ any of the sets of numbers $\{q_{\chi, P}\}_{\chi \in
G^\vee}$ which satisfy \eqref{eqn-red-cond-num}.
\end{enumerate}

\begin{proof} Take an open cover $U_i$ on which all 
$\mathcal{L}(-D_\chi)$ are trivialised and write $g_{\chi,i}$ for
the generator of $\mathcal{L}(-D_\chi)$ on $U_i$. $\{D_\chi\}$ being a
reductor set is equivalent to $\oplus \mathcal{L}(-D_\chi)$ being
closed under $\regring$-action on $K(\mathbb{C}^n)$. As $\regring$ is
a direct sum of its $G$-homogeneous parts, it is sufficient to check
the closure under the action of just the homogeneous functions. So on
each $U_i$, we want \begin{align*} f g_{\chi,i} \in \mathcal{O}_Y(U_i)
g_{\chi \rho(f),i} \end{align*} to hold for all $f \in \regring_G$,
$\chi \in G^\vee$. 

	On the other hand, with the notation above, $G$-Cartier
divisor $D_\chi + (f) - D_{\chi \rho(f)}$ is given on $U_i$ by
$\frac{f g_{\chi,i} }{g_{\chi\rho(f),i}}$ and it being effective is
equivalent to \begin{align*} \frac{f g_{\chi,i} }{g_{\chi \rho(f),i}}
\in \mathcal{O}_Y(U_i) \end{align*} for all $U_i$'s. 

	The result now follows.  \end{proof}

 We now translate the reductor condition \eqref{eqn-reductor-condition} 
into toric language and investigate what it implies for the reductor pieces 
of the family on the open toric charts $A_\sigma$ of a toric resolution $Y$.

	\begin{exmpl} \label{exmpl-canon-fam-red} Let $G$ and $Y$ be
as in previous examples. Let $\{D_\chi\}$ be a prereductor set where
each $D_\chi = \sum q_{\chi,i} E_i$  is given as follows
\begin{align*} 
& D_{\trch_}  =   0  &  & 
D_{\chi_1}  =   \frac{1}{8} E_4 + \frac{2}{8} E_5 + \frac{4}{8} E_6 + \frac{5}{8} E_7  & \\ 
& D_{\chi_2}  =    \frac{2}{8} E_4 + \frac{4}{8} E_5 + \frac{2}{8} E_7
&  
& D_{\chi_3}  =    \frac{3}{8} E_4 + \frac{6}{8} E_5 +
\frac{4}{8} E_6 + \frac{7}{8} E_7  & \\ 
& D_{\chi_4}  =  \frac{4}{8}
E_4 + \frac{4}{8} E_7  & 
& D_{\chi_5}  =  \frac{5}{8} E_4 + \frac{2}{8} E_5 + \frac{4}{8} E_6 +
\frac{1}{8} E_7  & \\
& D_{\chi_6} =  \frac{6}{8} E_4 + \frac{4}{8} E_5 + \frac{6}{8} E_7  & 
& D_{\chi_7}  =  \frac{7}{8} E_4 + \frac{6}{8} E_5 + \frac{4}{8} E_6 +
\frac{3}{8} E_7 & 
\end{align*}

In the view of Proposition \ref{prps-toric-valuation}, the reductor 
condition \eqref{eqn-reductor-condition} is equivalent to 
\begin{align} \label{prps-reductor-condition-toric}
q_{\chi,i} + e_i(m) - q_{\chi \rho(m), i} \geq 0 \end{align} for all
$\chi \in G^\vee$, $e_i \in \mathfrak{E}$ and $m \in \mathbb{Z}^n_+$. 

 The careful reader could now verify that
\eqref{prps-reductor-condition-toric} holds for $m = (1,0,0)$,
$(0,1,0)$ and $(0,0,1)$ and hence $\{D_\chi\}$ is a reductor set and
$\oplus\mathcal{L}( - D_\chi)$ is a family of $G$-constellations.

Recall now reductor pieces introduced in Definition \ref{subsection-reductors}. Let us calculate the reductor piece
$\{x^{p_\chi}\}$ specified by the generators of $\mathcal{L}(
- D_\chi)$ on the affine piece $A_{\left< e_5, e_6, e_7\right>}$. This 
is the same calculation of a generator of a $G$-Weil divisor on a
given open toric chart that we saw in Example \ref{exmpl-weil-divisor}, e.g.
$$ p_{\chi_1} = q_{\chi_1,5}\; \check{e}_5 + q_{\chi_1,6}\;
\check{e}_6 + q_{\chi_1,7}\; \check{e}_7$$ and so $$x^{p_{\chi_7}} =
\left( \frac{y^2 z}{x} \right)^{2/8} \left( \frac{z^2}{y}
\right)^{4/8} \left(\frac{x^2}{z^2} \right)^{5/8} = x $$ 

Repeating this for each $\chi \in G^\vee$, we obtain 
$\{ x^{p_\chi} \} = \{ 1, x, y, xy , \frac{x}{z}, z, \frac{xy}{z}, yz \}$, 
the reductor piece pictured below as a diagram in the monomial lattice
$\mathbb{Z}^n$:

\begin{center}
\includegraphics[scale=0.1]{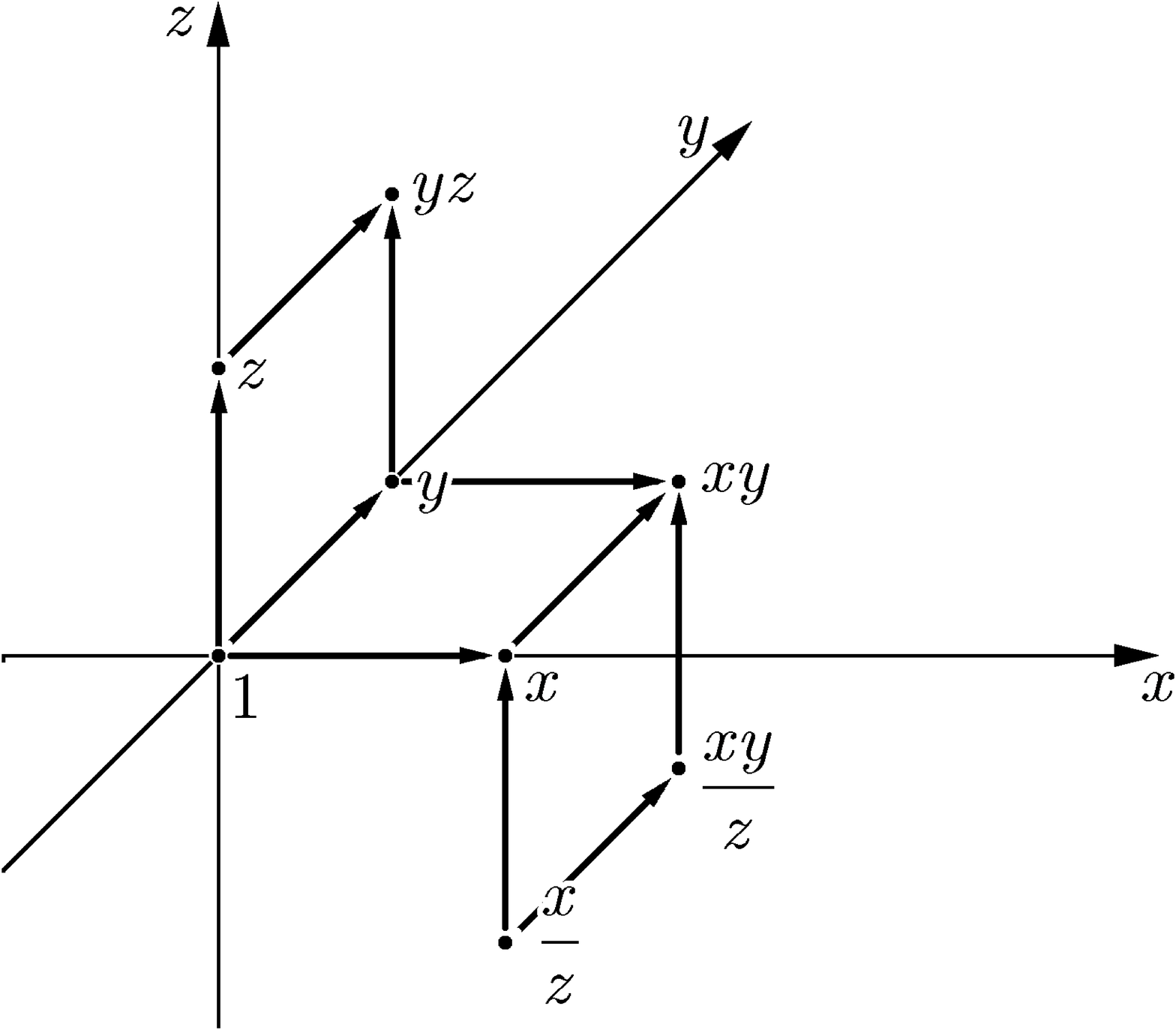}
\end{center}

The inequalities \eqref{prps-reductor-condition-toric} now translate into
the following form $$ e_i(p_\chi + m - p_{\chi\rho(m)}) > 0 \quad (i =
5,6,7)$$ that is \begin{align} \label{eqn-red-cond-piece}
\frac{x^{p_\chi} x^m}{x^{p_{\chi\rho(m)}}} \in \mathbb{C}[\sigma^\vee]
\end{align} for every $m \in \mathbb{Z}^n_+$. This agrees with the
discussion in Section \ref{subsection-reductors}, where it is
precisely the condition for $\oplus \mathcal{O}_{A_\sigma} x^{p_\chi}$
to be a family of $G$-constellations parametrised by $A_\sigma$. 

  The reader may find the diagrams set in the monomial lattice 
$\mathbb{Z}^n$ convenient for checking if a given monomial set 
$\{x^{p_\chi}\}$ satisfies the reductor equations in the form  
\eqref{eqn-red-cond-piece}. One merely needs to check that when adding 
$(1,0,0)$, $(0,1,0)$ or $(0,0,1)$ to any $p_\chi$, the vector reducing 
the result to $p_{\chi'}$ (for appropriate $\chi'$) lies within 
the cone $\sigma^\vee$. 

\begin{center} 
\includegraphics[scale=0.1]{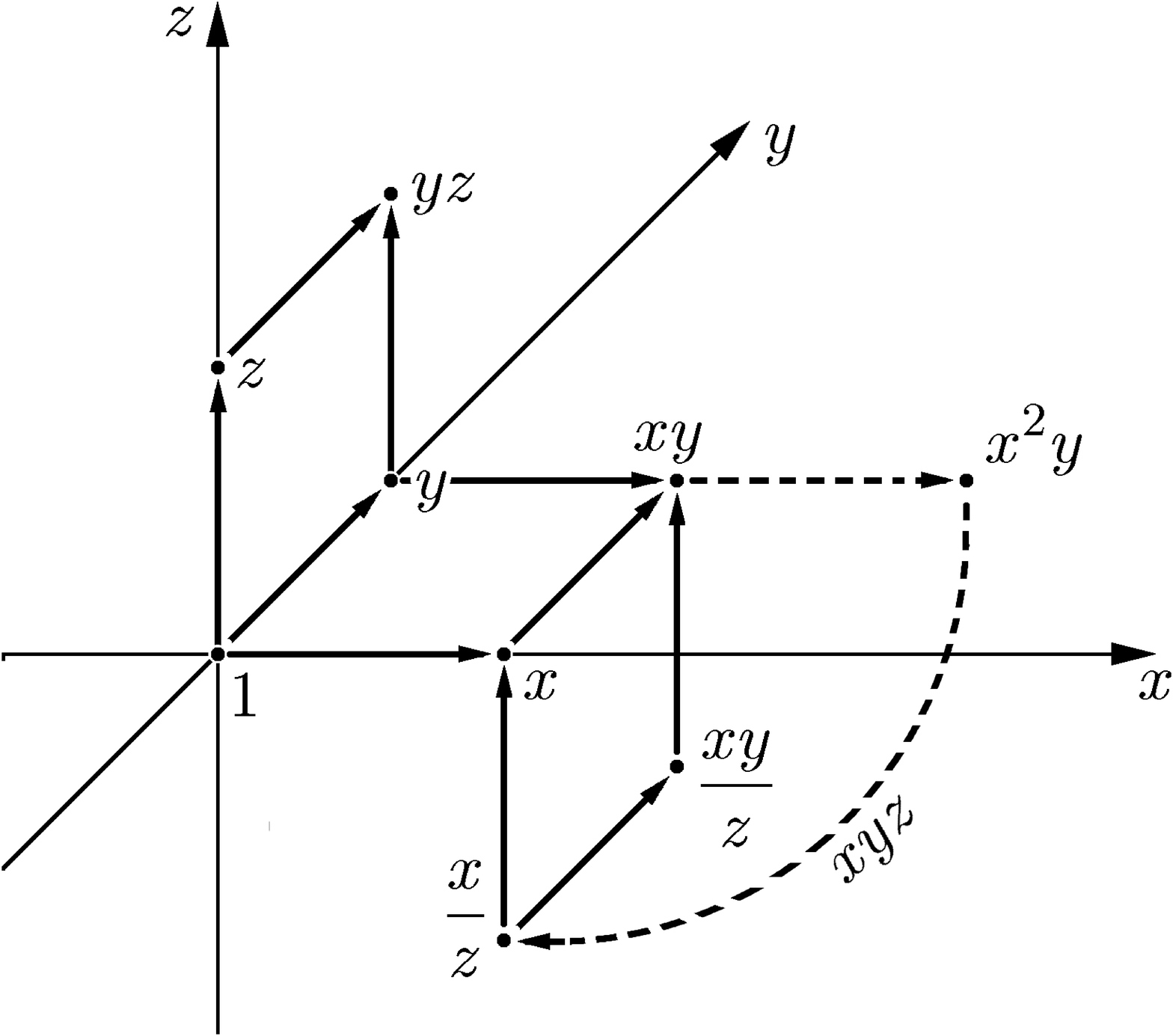}
\end{center}

\end{exmpl}

\subsection{Existence and symmetries} \label{section-canonical}

   So far we have seen no indication that, over an arbitrary
resolution $Y$, there exist any deformations of the generic orbit 
in the first place. There is apriori, for an arbitrary $Y$, 
no reason why it could at all be able to parametrise a family 
of $G$-constellations with such a strong relation to the geometry 
of $Y$ (Proposition \ref{prps-tfae}) as that of a deformation of 
the generic orbit. However, the following result shows that, for 
an absolutely any resolution $Y$, we always have at least one such family. 

\begin{prps}[Canonical family] 

    For an arbitrary resolution $Y$ of $X$ the set of $G$-Weil divisors
given by $D_\chi = \sum v(P, \chi) P$, where $P$ runs over all prime
Weil divisors on $Y$, satisfies the reductor condition.

    We shall call the family $\mathcal{F} = \oplus
\mathcal{L}(-D_\chi)$ \tt the canonical deformation of the generic
orbit of $G$ across $Y$\it.  

\end{prps} 

\bf Remark: \rm For $D_\chi = \sum v(P, \chi) P$ to be a $G$-Weil
divisor we need, in particular, for it to be a finite sum. This is
implied by Corollary \ref{prps-fract-parts-zero}.

\begin{proof} We need to show that for any
$\chi \in G^\vee$, any $G$-homogeneous $f \in \regring_G$ and any prime
divisor $P$ on $Y$ we have 

\begin{align*} 
v(P, \chi) + v_P(f) - v(P, \chi \rho(f)) \geq 0 
\end{align*}

First observe that the above expression must be
integer valued. Also $v(P, \chi)\geq 0$ and $ - v(P, \chi \rho(f))  >
-1 $ by definition, while $v_P(f) \geq 0$ since $f^n$ is regular on
all of $Y$. So we must have $$ v(P, \chi) + v_P(f) - v(P, \chi
\rho(f)) > -1 $$ and the result follows.  \end{proof}

\begin{cor} 

  Let $Y$ be a toric resolution of $X$. Then the canonical family of
$G$-constellations on $Y$ is given by $\{D_\chi\}$ where
\begin{align*} D_\chi = \sum_{i \in \mathfrak{E}} v(E_i, \chi) E_i
\end{align*}

   Moreover, on any affine open piece $A_\sigma$, we have
\begin{align} \label{canon-family-local-shape} 
\mathcal{F}(A_\sigma) = \mathbb{C}[\sigma \cap \mathbb{Z}^n] 
\end{align} 

\end{cor}

\begin{proof} 

  The first statement follows trivially from the definition of the
canonical family and the fact that $v(P,\chi) = 0$ whenever $P$ is not
one of the divisors $E_i$ (Corollary \ref{cor-except-nonzero}). 

   For the second statement, without loss of
generality let $\sigma = <e_1, \dots, e_n>$. Write $\mathcal{F}(A_\sigma)
= \oplus \mathbb{C}[\sigma^\vee \cap M] x^{p_\chi}$, where
$x^{p_\chi}$ are the generators of $\mathcal{L}(-D_\chi)(A_{\sigma})$. 
Proposition \ref{prps-gens-of-weil} implies that for each $p_\chi$ 
we have $e_i(p_\chi) = v(E_i, \chi)$ for all $i \in 1, \dots, n$. 
   But all the numbers $v(E_i, \chi)$ are positive by 
definition, which implies that each $p_\chi$ lies in $\sigma^\vee$  
and so $\mathcal{F}(A_\sigma)
\subseteq \mathbb{C}[\sigma^\vee \cap \mathbb{Z}^n]$. Conversely, 
given any $m \in \sigma^\vee \cap \mathbb{Z}^n$
$$e_i(m - p_{\rho(m)}) = e_i(m) - v(\rho(m), E_i) \geq 0$$ as 
$v(E_i), \rho(m)$ is precisely  the fractional part of 
$v_{E_i}(m) = e_i(m)$. Therefore $m  - p_{\rho(m)} \in \sigma^\vee \cap
M$
and so we have the inclusion in the other direction.  

Geometrically, one could easily convince oneself in the truth 
of this statement by picturing the cone 
$\sigma^\vee = \{ v \in \mathbb{R}^n \; | \; e_i(v) \geq 0 \}$ in 
$\mathbb{Z}^n \otimes \mathbb{R}$ 
and observing that the set $\{p_\chi\}$ of the exponents 
of the reductor piece of $\mathcal{F}$ on $A_\sigma$ consists precisely 
of all the elements of $\mathbb{Z}^n$ lying within the topmost area
$U$ of $\sigma^\vee$ given by $1 > e_i(v) \geq 0$. $\sigma^\vee
\cap \mathbb{Z}^n$ is then precisely $(U \cap \mathbb{Z}^n) + 
(\sigma^\vee \cap M)$.  
We can also see why reductor condition holds: as the cone 
$\mathbb{R}^n_+$ lies within the cone $\sigma^\vee$, 
$p_\chi + m$ lies within $\sigma^\vee \cap \mathbb{Z}^n$ for any 
$x^m \in \regring$. 

\end{proof}

\begin{exmpl}
The reductor set $\{D_\chi\}$ given in Example \ref{exmpl-canon-fam-red} 
specifies the canonical family on $Y$. Indeed, observe that all 
the numbers $q_{\chi,i}$ are between $0$ and $1$. The 
\eqref{eqn-gweil-cond} in definition of a $G$-Weil divisor implies 
they must be $v(E_i,\chi)$. 
  
Generally, to calculate the canonical family in a toric case, one needs 
to choose a monomial $m_\chi$ of weight $\chi$ for each $\chi \in G$. 
Then, for each $e_i \in \mathfrak{F}$, one calculates
the rational number $e_i(m_\chi)$ and takes its fractional part, 
which is precisely $v(E_i, \chi)$. The G-Weil
divisors $D_\chi = \sum_i  v(E_i, \chi) E_i$ are then the reductor
set for the canonical family.

For instance, the numbers for the canonical family in 
Example \ref{exmpl-canon-fam-red} were obtained as follows: take character 
$\chi_3 \in G^\vee$ and then take $x^3$, a monomial of weight $\chi_3$. 
Calculating $e_5(3,0,0) = \frac{1}{8}(2 * 3 +
4 * 0 + 2 * 0) = \frac{6}{8}$, we obtain the coefficient of $E_5$ in
$D_{\chi_3}$. Similarly $e_7(3,0,0) = \frac{15}{8}$ and its fractional
part $\frac{7}{8}$ is the coefficient of $E_7$ in $D_{\chi_3}$. 

Observe also that given any other reductor set $\{D_\chi'\}$, its 
$q_{\chi',i}$ will differ from those of the canonical one by integer
numbers. 

Observe also that on the level of reductors $\{x^{p_\chi}\}$, 
the change introduced to the family by adding an integer $n$ to 
$q_{\chi,i}$ amounts precisely to shifting $p_\chi$ by $n \check{e}_i$ in 
the reductor of those open pieces $A_\sigma$ where $e_i \in \sigma$.  
But note that $\check{e}_i$ is a different vector in $M$ for each such 
$\sigma$.      
\end{exmpl}

Having established that deformations of the generic orbit across $Y$
always exist, we now consider symmetries which the set of them  
must possess.

\begin{prps}[Character Shift]
Let $\{D_\chi\}$ be a reductor set. Then for any $\lambda$-Weil divisor $N$, 
the set $\{ D_\chi + N \}$ also satisfies the reductor condition. 

Moreover, up to equivalence of families, the deformation 
$\mathcal{F}'$ it specifies depends only on $\lambda$ and not on the
choice of $N$, and the unique normalized reductor set $\{D'_\chi\}$ 
specifying $\mathcal{F}'$ is given by 
\begin{align}
D'_{\chi  \lambda} = D_{\chi} - D_{\lambda^{-1}}
\end{align}
\end{prps}
\begin{proof}
That the new set of divisors satisfies the reductor condition 
is trivial:
$$ (D_\chi + N) + (m) - (D_{\chi \rho(m)} + N) \geq 0 $$
is immediately equivalent to the statement that $\{D_\chi\}$ satisfy 
the reductor condition.

 For the second claim, observe that the divisor in the trivial 
character class is now $(D_{\lambda^{-1}} + N)$. Normalising by it 
we obtain in character class $\chi + \lambda$ 
$$ D_{\chi} + N - D_{\lambda^{-1}} - N $$
which establishes the claim.
\end{proof}

\begin{defn}
Given a normalized reductor set $\{D_\chi\}$, we shall call 
normalized reductor set $\{D_{\chi} - D_{\lambda^{-1}}\}$ 
the \tt $\lambda$-shift \rm of $\{D_\chi\}$. 
\end{defn}

\begin{exmpl}
On the level of reductors $\{x^{p_\chi}\}$, $\lambda$-shift leaves 
the geometrical configuration of $p_\chi$'s in the lattice $\mathbb{Z}^n$ 
the same, but permutes them and shifts the origin to the
new location of $p_\trch_$.  

For example, consider the case of the reductor piece calculated in Example 
\ref{exmpl-canon-fam-red}. After a $\chi_4$-shift it becomes:

\begin{center}
\includegraphics[scale=0.1]{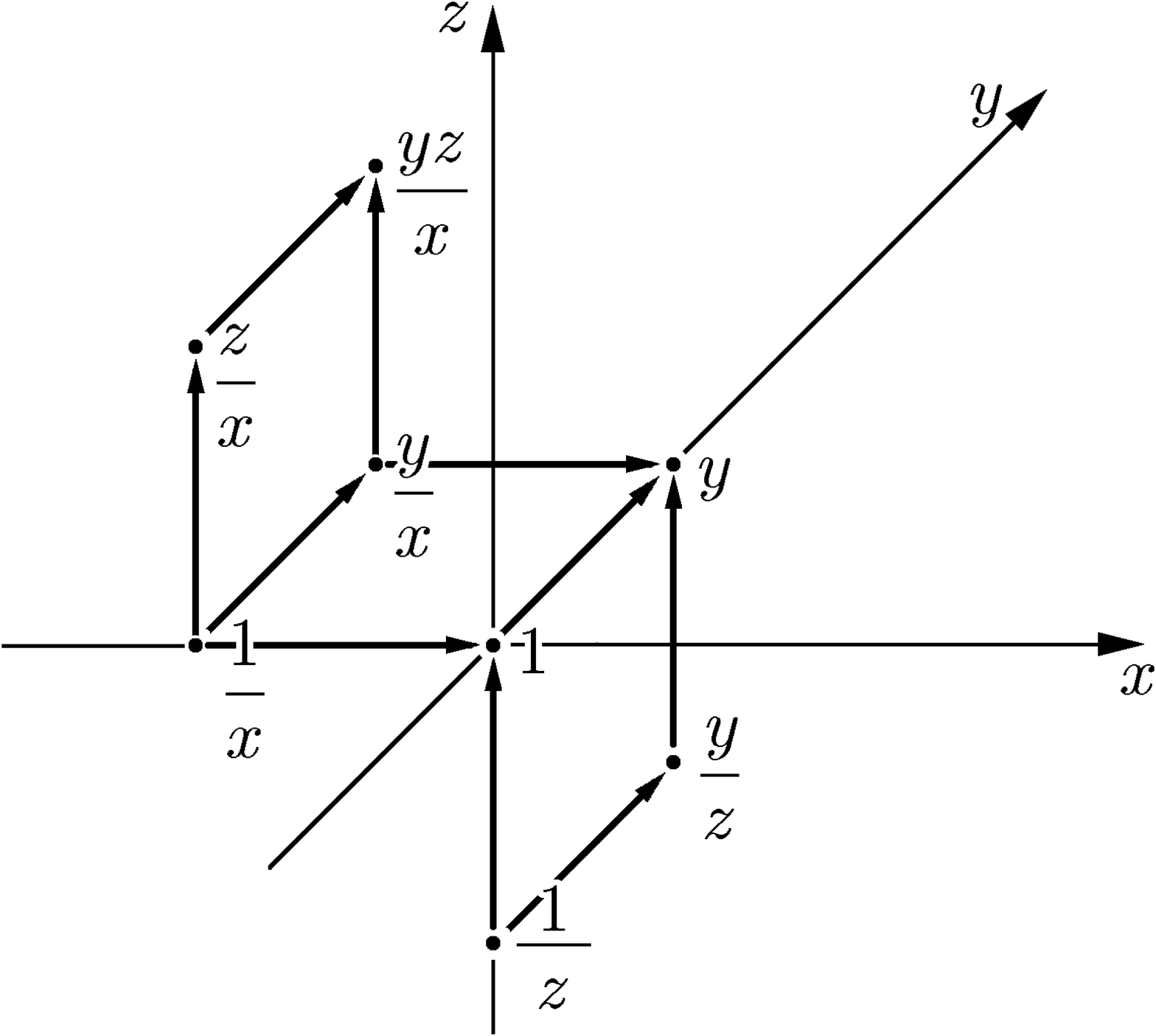}
\end{center}

\end{exmpl}

\begin{prps}[Reflection] 
Let $\{D_\chi\}$ be a reductor set. 
Then the set $\{ - D_\chi \}$ also satisfies the reductor condition.  
\end{prps}
\begin{proof}
We need to show that 
\begin{align*}
- D_{\chi^{-1}} + (m) - ( - D_{\chi^{-1}\rho(m)^{-1}}) \geq 0
\end{align*}
Rearranging we get
\begin{align*}
D_{\chi^{-1} \rho(m)^{-1}} + (m) - D_{\chi^{-1} \rho(m)^{-1} \rho(m)}
\geq 0
\end{align*}
which is one of the reductor equations the original set $\{D_\chi\}$
must satisfy.
\end{proof}

\begin{defn}
Given a reductor set $\{D_\chi\}$, we shall call the reductor set 
$\{-D_{\chi}\}$ the \tt reflection \rm of $\{D_\chi\}$.
\end{defn}

\begin{exmpl}

On the level of reductors $\{x^{p_\chi}\}$, the reflection is precisely 
the reflection of $p_\chi$ about the origin in the lattice $\mathbb{Z}^n$.

For example, consider the case of the reductor piece calculated in 
Example \ref{exmpl-canon-fam-red}. After a reflection it becomes:

\begin{center}
\includegraphics[scale=0.1]{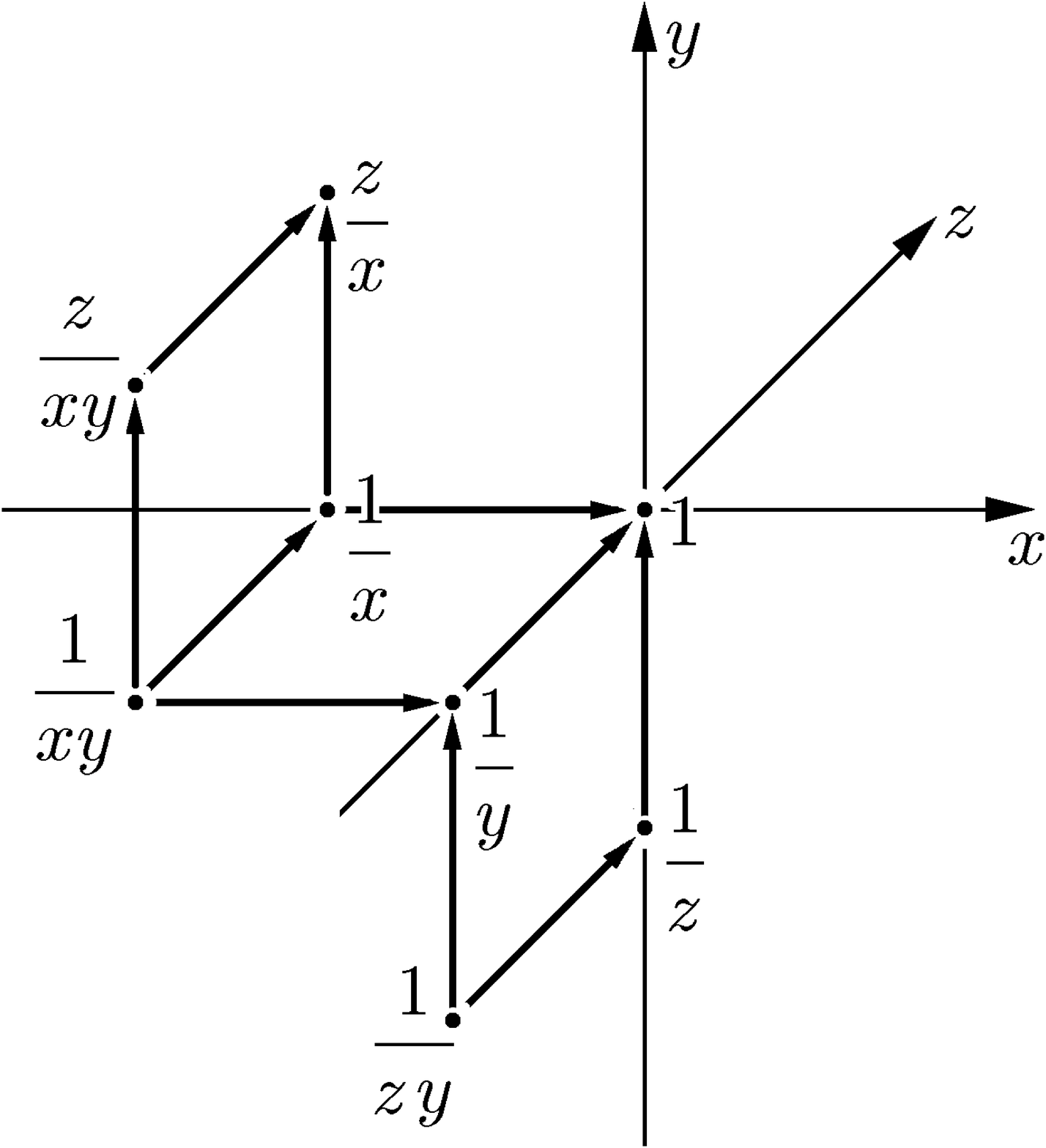}
\end{center}

\end{exmpl}

\subsection{Maximal Shifts}

  We now examine the individual line bundles $\mathcal{L}(-D_\chi)$ in 
a deformation of the generic orbit and show that the reductor condition 
imposes a restriction on how far apart from each other they can be. 

\begin{lemma} \label{lemma-gen-max-shift}
 Let $\{D_\chi\}$ be a reductor set. Write each $D_\chi$ as 
$\sum q_{\chi, P} P$, where $P$ ranges over all the prime Weil
divisors on $Y$. Then we necessarily have for any 
$\chi_1, \chi_2 \in G^\vee$ and for any prime Weil divisor $P$
\begin{align}
\min_{f \in \regring_{\chi_1/\chi_2}} v_P(f),   
\quad \geq \quad q_{\chi_1,P} - q_{\chi_2,P}
\quad \geq \quad 
- \min_{f \in \regring_{{\chi_2}/{\chi_1}}} v_P(f)
\end{align}
where $\regring_\chi$ is the set of all the $\chi$-homogeneous 
functions in $\regring$. 
\end{lemma}
\begin{proof}
Both inequalities follow directly from the reductor condition 
\eqref{eqn-reductor-condition}: the right inequality by setting 
$\chi = \chi_1 \in G^\vee$, $\rho(f) = \frac{\chi_2}{\chi_1}$ 
and letting $f$ vary within $\regring_{\rho(f)}$; the left inequality
by setting $\chi = \chi_2$ and $\rho(f) = \frac{\chi_1}{\chi_2}$.
\end{proof}

This suggests the following definition:
\begin{defn} \label{def-maxshift}
For each character $\chi \in G^\vee$, the 
\tt maximal shift $\chi$-divisor $M_\chi$ \it is defined to be
\begin{align} \label{eqn-max-shift-defn}
M_\chi = \sum_P (\min_{f \in \regring_\chi} v_P(f)) P
\end{align}
where $P$ ranges over all prime Weil divisors on $Y$.  
\end{defn}

Observe that the fact that the sum in \eqref{eqn-max-shift-defn} is finite 
follows directly from Corollary \ref{cor-except-nonzero}.

\begin{lemma} \label{lemma-max-shift-is-reductor}
The $G$-Weil divisor set $\{M_\chi\}$ is a normalised reductor set. 
\end{lemma}
\begin{proof}
To show that the set $ \{ M_\chi \} $ satisfies the reductor condition,  
we need to show that for every $f \in \regring_G$ 
and any prime divisor $P$ on $Y$
$$ v_P(m_\chi) + v_P(f) - v_P(m_{\chi \rho(f)}) \geq 0 $$
where $m_\chi$ and $m_{\chi \rho(f)}$ are chosen to achieve the 
minimality in \eqref{eqn-max-shift-defn}. 

 Observe that $m_\chi f$ is also a $G$-homogeneous element of $\regring$,
therefore by the minimality of $v_P (m_{\chi \rho(f)})$ we have
$$ v_P(m_\chi f)  \geq  v_P (m_{\chi \rho(f)}) $$
as required. 

To establish that $M_\trch_ = 0$, we observe that $v_P(1) = 0$ for
any prime Weil divisor $P$ on $Y$ and $v_P(f) \geq 0$ for any 
$G$-homogeneous $f \in \regring$. 

\end{proof}

Observe that with Lemma \ref{lemma-max-shift-is-reductor} we have
established another deformation of the generic orbit of $G$ which always
exists across any resolution $Y$. While in some cases it coincides 
with the canonical family, the reader will see in Example 
\ref{exmpl-maxshift-calc} the case when the canonical family and the maximal
shift family differ. 

Putting together Lemmas \ref{lemma-gen-max-shift} and 
\ref{lemma-max-shift-is-reductor} gives a result which 
shows that that the reductor set $\{M_\chi\}$ and its reflection $\{-M_\chi\}$ 
provide bounds on the set of all normalized reductor sets on $Y$.  

\begin{prps}[Maximal Shifts] \label{prps-max-shifts}
Let $\{D_\chi\}$ be a normalized reductor set.
Then for any $\chi \in G^\vee$  
\begin{align} \label{eqn-family-bounded}
M_\chi \geq
D_\chi \geq
- M_{\chi^{-1}}
\end{align}
Moreover both the bounds are achieved.
\end{prps}
\begin{proof}
To establish that \eqref{eqn-family-bounded} holds, set $\chi_2
= \trch_$ in Lemma \ref{lemma-gen-max-shift}. Lemma 
\ref{lemma-max-shift-is-reductor} shows that bounds are achieved.
\end{proof}

\begin{exmpl} \label{exmpl-maxshift-calc}
Let us calculate the maximal shift divisor set $\{M_\chi\}$ for 
the setup introduced in the Example \ref{exmpl-the-setup}. 

  By the definition $M_\chi = \sum m_{\chi,P} P  $ where $m_{\chi,P}
= \min_{f \in \regring_\chi} v_P(f)$. By Corollary 
\ref{cor-except-nonzero}, the numbers $m_{\chi,P}$ are only non-zero 
for divisors corresponding to elements of $\mathfrak{E}$. Therefore for 
each $e_i \in \mathfrak{E}$, we need to find $m_{\chi, E_i} = \min e_i(p)$ 
where $p$ ranges over elements of $\mathbb{Z}^n_+$ such that $\rho(p) = \chi$.

  It is only necessary to consider a finite number of choices for $p$ to
establish each $m_{\chi,P}$. Observe that it suffices to take the ones with 
$0 \leq p_i \leq |G|$, as $p' = p - (0,\dots, 0, |G|,0 , \dots, 0)$ is
again element of $\mathbb{Z}^n$ with $\rho(p') = \rho(p)$ and $e_i(p')
\leq e_i(p)$ for all $e_i \in \mathfrak{E}$. 

  For example, taking $e_5 = \frac{1}{8}(2,4,2)$ and considering all 
such $p$ we see that:
\begin{align*}
&m_{\trch_, E_5} = v_{E_5}(1) = e_5(0,0,0) = 0 
&m_{\chi_1, E_5} = v_{E_5}(x) = e_5(1,0,0) = \frac{2}{8}  \\
&m_{\chi_2, E_5} = v_{E_5}(x^2) = e_5(2,0,0) = \frac{4}{8}  
&m_{\chi_3, E_5} = v_{E_5}(x^3) = e_5(3,0,0) = \frac{6}{8} \\
&m_{\chi_4, E_5} = v_{E_5}(x^4) = e_5(4,0,0) = 1 
&m_{\chi_5, E_5} = v_{E_5}(z) = e_5(0,0,1) = \frac{2}{8} \\
&m_{\chi_6, E_5} = v_{E_5}(zx) = e_5(1,0,1) = \frac{4}{8}
&m_{\chi_7, E_5} = v_{E_5}(zx^2) = e_5(2,0,1) = \frac{6}{8}
\end{align*}
Observe that in case of $\chi_4$ we have $m_{P,\chi} \neq v_{P,\chi}$. 
So the maximal shift family for this $Y$ differs from the canonical family.  

  If we repeat this calculation for all elements of $\mathfrak{E}$, 
to obtain all numbers $m_{e_i,\chi}$, we will obtain:
\begin{align*}
&M_{\trch_} = 0, &M_{\chi_1} = \frac{1}{8} E_4 + \frac{2}{8} E_5 + \frac{4}{8} E_6 +\frac{5}{8} E_7 \\ 
&M_{\chi_2} = \frac{2}{8} E_4 + \frac{4}{8} E_5 +\frac{2}{8} E_7 &M_{\chi_3} = \frac{3}{8} E_4 + \frac{6}{8} E_5 + \frac{4}{8} E_6
+\frac{7}{8} E_7 \\ 
&M_{\chi_4} = \frac{4}{8} E_4 + E_5 + \frac{4}{8} E_7 &M_{\chi_5} = \frac{5}{8} E_4 + \frac{2}{8} E_5 + \frac{4}{8} E_6
+\frac{1}{8} E_7 \\ 
&M_{\chi_6} = \frac{6}{8} E_4 + \frac{4}{8} E_5 +\frac{6}{8} E_7 &M_{\chi_7} = \frac{7}{8} E_4 + \frac{6}{8} E_5 + \frac{4}{8} E_6
+\frac{3}{8} E_7 
\end{align*}
  Compare it to the reductor set of the canonical family given in 
Example \ref{exmpl-canon-fam-red}. 

  If we now wanted to calculate all the normalised reductor sets
(and hence all the normalised deformations of the generic orbit), we
simply need to check each of the finite number of 
prereductor sets between $\{M_\chi\}$ and its reflection
$\{-M_\chi\}$ and pick out the ones which satisfy the reductor
condition \eqref{eqn-reductor-condition}. 

  Recall now the remark after Proposition \ref{prps-reductor-condition}, 
about checking reductor condition independently at each prime divisor 
in $Y$. Here, it means that for any reductor set 
$\{\sum_i q_{\chi,i} E_i\}_{\chi \in G^\vee}$, the numbers 
$\{q_{\chi,i}\}_{\chi \in G^\vee}$ satisfy or fail the reductor
condition inequalities independently for each $e_i \in \mathfrak{E}$. 
This can be seen from the fact that each of the inequalities 
$\eqref{prps-reductor-condition-toric}$ features numbers 
$q_{\chi,i}$ all for the same $i$. 

  In particular it means that to list all the possible normalized
reductor sets on $Y$, it is sufficient to list for each 
$E_i$ all the sets $\{q_{\chi,i}\}_{\chi \in G^\vee}$
satisfying the inequalities $\eqref{prps-reductor-condition-toric}$. 
Then all the normalized reductor sets on $Y$ are given by all
the possible choices of one of these sets $\{q_{\chi,i}\}_{\chi \in G^\vee}$
for each $E_i$.

  For our particular $Y$, we give such list below:  

$E_4:$ \\
\begin{align*}
\frac{1}{8} \left( \begin{tabular}{r r r r r r r r}
$\trch_ $ & $ \chi_1 $ & $ \chi_2 $ & $ \chi_3 $ & $ \chi_4 $
& $ \chi_5 $ & $ \chi_6 $ & $ \chi_7 $ \\
 $0$ & $ 1 $ & $ 2 $ & $ 3 $ & $
4 $ & $ 5 $ & $ 6 $ & $ 7 $ \\
 $0$ & $ 1 $ & $ 2 $ & $ 3 $ & $
4 $ & $ 5 $ & $ 6 $ & $ - 1 $ \\
 $0$ & $ 1 $ & $ 2 $ & $ 3 $ & $
4 $ & $ 5 $ & $ - 2 $ & $ - 1 $ \\
 $0$ & $ 1 $ & $ 2 $ & $ 3 $ & $
4 $ & $ - 3 $ & $ - 2 $ & $ - 1 $ \\
 $0$ & $ 1 $ & $ 2 $ & $ 3 $ & $
- 4 $ & $ - 3 $ & $ - 2 $ & $ - 1 $ \\
 $0$ & $ 1 $ & $ 2 $ & $ - 5 $ & $
- 4 $ & $ - 3 $ & $ - 2 $ & $ - 1 $ \\
 $0$ & $ 1 $ & $ - 6 $ & $ - 5 $ & $
- 4 $ & $ - 3 $ & $ - 2 $ & $ - 1 $ \\
 $0$ & $ - 7 $ & $ - 6 $ & $ - 5 $ & $
- 4 $ & $ - 3 $ & $ - 2 $ & $ - 1 $ \\
\end{tabular} \right)
\end{align*} 

$E_5:$ \\
\begin{align*}
\frac{1}{8}
\left( \begin{tabular}{r r r r r r r r}
$\trch_ $ & $ \chi_1 $ & $ \chi_2 $ & $ \chi_3 $ & $ \chi_4 $
& $ \chi_5 $ & $ \chi_6 $ & $ \chi_7 $ \\
 $0$ & $ 2 $ & $ 4 $ & $ 6 $ & $ 8 $ & $ 2 $ & $ 4 $ & $ 6 $ \\
 $0$ & $ 2 $ & $ 4 $ & $ 6 $ & $ 0 $ & $ 2 $ & $ 4 $ & $ 6 $ \\
 $0$ & $ 2 $ & $ 4 $ & $ - 2 $ & $ 0 $ & $ 2 $ & $ 4 $ & $ 6 $ \\
 $0$ & $ 2 $ & $ 4 $ & $ 6 $ & $ 0 $ & $ 2 $ & $ 4 $ & $ - 2 $ \\
 $0$ & $ 2 $ & $ 4 $ & $ - 2 $ & $ 0 $ & $ 2 $ & $ 4 $ & $ - 2 $ \\
 $0$ & $ 2 $ & $ - 4 $ & $ - 2 $ & $ 0 $ & $ 2 $ & $ 4 $ & $ - 2 $ \\
 $0$ & $ 2 $ & $ 4 $ & $ - 2 $ & $ 0 $ & $ 2 $ & $ - 4 $ & $ - 2 $ \\
 $0$ & $ 2 $ & $ - 4 $ & $ - 2 $ & $ 0 $ & $ 2 $ & $ - 4 $ & $ - 2 $ \\
 $0$ & $ - 6 $ & $ - 4 $ & $ - 2 $ & $ 0 $ & $ 2 $ & $ - 4 $ & $ - 2 $ \\
 $0$ & $ 2 $ & $ - 4 $ & $ - 2 $ & $ 0 $ & $ - 6 $ & $ - 4 $ & $ - 2 $ \\
 $0$ & $ - 6 $ & $ - 4 $ & $ - 2 $ & $ 0 $ & $ - 6 $ & $ - 4 $ & $ - 2 $ \\
 $0$ & $ - 6 $ & $ - 4 $ & $ - 2 $ & $ - 8 $ & $ - 6 $ & $ - 4 $ & $ - 2 $ \\
\end{tabular} \right)
\end{align*}

$E_6:$ \\
\begin{align*}
\frac{1}{8}
\left( 
\begin{tabular}{r r r r r r r r}
$\trch_ $ & $ \chi_1 $ & $ \chi_2 $ & $ \chi_3 $ & 
$ \chi_4 $ & $ \chi_5 $ & $ \chi_6 $ & $ \chi_7 $ \\
$0 $ & $ 4 $ & $ 0 $ & $ 4 $ & $ 0 $ & $ 4 $ & $ 0 $ & $ 4$ \\
$0 $ & $ -4 $ & $ 0 $ & $ -4 $ & $ 0 $ & $ -4 $ & $ 0 $ & $ -4$ 
\end{tabular}
\right)
\end{align*}

$E_7:$ \\
\begin{align*}
\frac{1}{8}
\left( 
\begin{tabular}{r r r r r r r r}
$\trch_ $ & $ \chi_1 $ & $ \chi_2 $ & $ \chi_3 $ & $ \chi_4 $ & $
\chi_5 $ & $ \chi_6 $ & $ \chi_7 $ \\
$0 $ & $ 5 $ & $ 2 $ & $ 7 $ & $ 4 $ & $ 1 $ & $ 6 $ & $ 3$\\
$0 $ & $ 5 $ & $ 2 $ & $ -1 $ & $ 4 $ & $ 1 $ & $ 6 $ & $ 3$\\
$0 $ & $ 5 $ & $ 2 $ & $ -1 $ & $ 4 $ & $ 1 $ & $ -2 $ & $ 3$\\
$0 $ & $ -3 $ & $ 2 $ & $ -1 $ & $ 4 $ & $ 1 $ & $ -2 $ & $ 3$\\
$0 $ & $ -3 $ & $ 2 $ & $ -1 $ & $ -4 $ & $ 1 $ & $ -2 $ & $ 3$\\
$0 $ & $ -3 $ & $ 2 $ & $ -1 $ & $ -4 $ & $ 1 $ & $ -2 $ & $ -5$\\
$0 $ & $ -3 $ & $ -6 $ & $ -1 $ & $ -4 $ & $ -7 $ & $ -2 $ & $ -5$
\end{tabular}
\right)
\end{align*} 
\end{exmpl}

For one particular resolution $Y$, the family provided by 
the maximal shift divisors is already quite well-known. 

\begin{prps} \label{prps-maxshift-ghilb}
Let $Y = G$-$\hilb \mathbb{C}^n$, the moduli space of
$G$-clusters in $\mathbb{C}^n$. If $Y$ is smooth, then 
$\oplus \mathcal{L}(-M_\chi)$ is the universal 
family $\mathcal{F}$ of $G$-clusters parametrised by
$Y$, up to the usual equivalence of families. 
\end{prps}

\begin{proof}
Firstly $\mathcal{F}$ is a deformation of the generic orbit, as
over any set $U \subset X$ such that $G$ acts freely on $q^{-1}(U)$
we have $\pi_* \mathcal{F} |_U \simeq q_* \mathcal{O}_{\mathbb{C}^n}
|_U$.
Hence write $\mathcal{F}$ as $\oplus \mathcal{L}(-D_\chi)$ for some
reductor set $\{ D_\chi \}$. Take an open cover $\{U_i\}$ of $Y$ and
consider the generators $\{f_{\chi,i}\}$ of $D_\chi$ on each
$U_i$. Working up to equivalence, we can consider $\{D_\chi\}$ 
to be normalised and so $f_{\trch_, i} = 1$ for all $U_i$.

Now any $G$-cluster $Z$ is given by some invariant ideal $I \subset
\regring$ and so the corresponding $G$-constellation
$H^0(\mathcal{O}_Z)$
is given by $\regring / I$. In particular note that $R/I$ is generated
by $R$-action on the generator of $\trch_$-eigenspace. Therefore
any $f_{\chi, i}$ is generated from $f_{\trch_, i} = 1$ by
$\regring$-action, which means that all $f_{\chi, i}$ lie in
$\regring$.

But this means that for any prime Weil divisor $P$ on $Y$ we have
\begin{align*}
v_P(f_{\chi, i}) \geq \min_{f \in \regring_\chi} v_P(f)
\end{align*}
and therefore $D_\chi \geq M_\chi$. Now Corollary \ref{prps-max-shifts}
forces the equality. 
\end{proof}

\subsection{Summary}

Finally, we combine the results achieved thus far 
into a classification theorem.

\begin{theorem}[Classification]
\label{theorem-classification}

Let $G$ be a finite abelian subgroup of $\gl_n(\mathbb{C})$, $X$ be the
quotient of $\mathbb{C}^n$ by the action of $G$ and $Y$ be a resolution
of $X$. Then all deformations of the generic orbit across $Y$, 
up to isomorphism, are of form 
$\oplus_{\chi \in G^\vee} \mathcal{L}(-D_\chi)$, where each 
$D_\chi$ is a $\chi$-Weil divisor and the set $\{D_\chi\}$ satisfies 
the inequalities:
\begin{align*}
D_\chi + (f) - D_{\chi \rho(f)} \geq 0 
\end{align*}
for all $\chi \in G^\vee$ and all $G$-homogeneous $f \in \regring$.
Here $\rho(f)$ is the homogeneous weight of $f$. Conversely for any 
such set $\{D_\chi\}$, $\oplus \mathcal{L}(-D_\chi)$ is a deformation
of the generic orbit.     

Moreover, each equivalence class of families has precisely one family 
with $D_{\chi_0} = 0$. The divisor set $\{D_\chi\}$ corresponding 
to such a family satisfies inequalities 
\begin{align*} M_\chi \geq
D_\chi \geq
- M_{\chi^{-1}}
\end{align*}
where $\{ M_{\chi} \}$ is a fixed divisor set depending only on $G$
and $Y$. In particular, the number of equivalence classes of families 
is finite. 
\end{theorem}

\begin{proof}
Proposition \ref{prps-reductor-condition} establishes 
the correspondence of isomorphism classes of deformations of 
the generic orbit and reductor sets. 
Corollary \ref{cor-equiv-classes-normalized-sets} lifts the
correspondence to the level of equivalence classes and normalised
reductor sets. Corollary \ref{prps-max-shifts} gives the bounds 
on the set of all normalised
reductor sets, and as due to Corollary \ref{cor-except-nonzero}
each $M_\chi$ is a finite sum, this set is finite. 
\end{proof}

\bibliography{references}

\providecommand{\bysame}{\leavevmode\hbox to3em{\hrulefill}\thinspace}
\providecommand{\MR}{\relax\ifhmode\unskip\space\fi MR }
\providecommand{\MRhref}[2]{%
  \href{http://www.ams.org/mathscinet-getitem?mr=#1}{#2}
}
\providecommand{\href}[2]{#2}
\begin{thebibliography}{BKR01}

\bibitem[Ben94]{Bens94}
D.~J. Benson, \emph{Polynomial invariants of finite groups}, Cambridge
  University Press, 1994.

\bibitem[BKR01]{BKR01}
T.~Bridgeland, A.~King, and M.~Reid, \emph{The {M}c{K}ay correspondence as an
  equivalence of derived categories}, J. Amer. Math. Soc. \textbf{14} (2001),
  535--554.

\bibitem[CI02]{Craw-Ishii-02}
A.~Craw and A.~Ishii, \emph{Flops of ${G}-\hilb$ and equivalences of derived
  category by variation of {GIT} quotient}, preprint math.AG/0211360, (2002).

\bibitem[CR02]{Craw02}
A.~Craw and M.~Reid, \emph{How to calculate ${A}-\hilb \mathbb{C}^3$},
  Seminaires et Congres \textbf{6} (2002), 129--154.

\bibitem[Cra01]{Craw-thesis}
A.~Craw, \emph{The {M}c{K}ay correspondence and representations of the
  {M}c{K}ay quiver}, Ph.D. thesis, University of {W}arwick, 2001.

\bibitem[Dan78]{Dani78}
V.I. Danilov, \emph{The geometry of toric varieties}, Russian Math. Surveys
  \textbf{33} (1978), 97--154.

\bibitem[Gar86]{Garl86}
D.J.H. Garling, \emph{A course in {G}alois theory}, Cambridge University Press,
  1986.

\bibitem[Har77]{Harts77}
R.~Hartshorne, \emph{Algebraic geometry}, Springer-Verlag, 1977.

\bibitem[IR96]{ItReid96}
Y.~Ito and M.~Reid, \emph{The {M}c{K}ay correspondence for the finite subgroups
  of $\gsl(3,\mathbb{C})$.}, Higher-dimensional complex varieties ({T}rento
  1994), de Gruyter, 1996, pp.~221--240.

\bibitem[Kin94]{King94}
A.~King, \emph{Moduli of representations of finite-dimensional algebras},
  Quart. J. Math. Oxford \textbf{45} (1994), 515--530.

\bibitem[Kro89]{Kron89}
P.~Kronheimer, \emph{The construction of {ALE} spaces as hyper-{K}ahler
  quotients}, J. Diff. Geom. \textbf{29} (1989), 665--683.

\bibitem[Mat86]{Mats86}
H.~Matsumura, \emph{Commutative ring theory}, Cambridge University Press, 1986.

\bibitem[Nak00]{Nak00}
I.~Nakamura, \emph{Hilbert schemes of abelian group orbits}, J. Alg. Geom.
  \textbf{10} (2000), 775--779.

\bibitem[Rei87]{YPG87}
Miles Reid, \emph{Young person's guide to canonical singularities}, Proc. of
  Symposia in Pure Math. \textbf{46} (1987), 345--414.

\bibitem[SI96a]{Sar-In96a}
A.~Sardo-Infirri, \emph{Partial resolutions of orbifold singularities via
  moduli spaces of {HYM}-type bundles}, preprint math.AG/9610004, (1996).

\bibitem[SI96b]{Sar-In96b}
\bysame, \emph{Resolutions of orbifold singularities and the transportation
  problem on the {M}c{K}ay quiver}, preprint math.AG/9610005, (1996).

\end{thebibliography}
\bibliographystyle{amsalpha}
\end{document}